\documentclass[a4paper,11pt,onecolumn,twoside]{article}
\usepackage{mathrsfs}
\usepackage{mathrsfs}
\usepackage{amsfonts}
\usepackage{fancyhdr}
\usepackage{amsmath,amsfonts,amssymb,graphicx,amsthm}
\allowdisplaybreaks

\usepackage{amsmath,amsfonts,amssymb,amscd,amsthm,amsbsy,bm,latexsym}
\usepackage{gensymb}
\usepackage[pdftex,bookmarksnumbered=true,bookmarksopen=true,          
            colorlinks=true,pdfborder=001,citecolor=blue,
            linkcolor=red,anchorcolor=green,urlcolor=blue]{hyperref}
\usepackage{mathrsfs}
\allowdisplaybreaks

\usepackage{color}

\begin{document}

\title{\bf Asymptotic stability of viscous shock  profiles  for  the 1D compressible Navier-Stokes-Korteweg system with boundary effect }
\author{
{\bf Zhengzheng Chen}\\
School of Mathematical Sciences,
Anhui University, Hefei, 230601,  China\\[2mm]
{\bf  Yeping Li\thanks{Corresponding author. E-mail:
yplee@ecust.edu.cn}}\\
Department of Mathematics, East China University of Science and Technology, Shanghai, China\\[2mm]
{\bf Mengdi Sheng}\\
School of Mathematical Sciences,
Anhui University, Hefei, 230601,  China\\[2mm]
}

\date{}

\vskip 0.2cm

\maketitle

\vskip 0.2cm \arraycolsep1.5pt
\newtheorem{Lemma}{Lemma}[section]
\newtheorem{Theorem}{Theorem}[section]
\newtheorem{Definition}{Definition}[section]
\newtheorem{Proposition}{Proposition}[section]
\newtheorem{Remark}{Remark}[section]
\newtheorem{Corollary}{Corollary}[section]

\begin{abstract}

 This paper is concerned with the time-asymptotic behavior of strong solutions to an initial-boundary value problem of the  compressible Navier-Stokes-Korteweg system on the half line $\mathbb{R}^+$. The asymptotic profile of  the problem is shown to be a shifted viscous shock profile, which is suitably away from the boundary.  Moreover, we prove that  if the initial data around the shifted viscous shock profile and the strength of the shifted viscous shock profile are sufficiently small, then the problem has a unique global strong solution, which tends  to the shifted viscous shock profile as time goes to infinity. The analysis
 is based  on the elementary $L^2$-energy method and the key point is to deal with the boundary estimates.

\bigbreak
\noindent

\bigbreak

\end{abstract}

\section{Introduction }
\setcounter{equation}{0}
 In this paper, we consider the following  one-dimensional compressible Navier-Stokes-Korteweg system  on the half line $\mathbb{R}^+=(0,+\infty)$ in the  Lagrangian  coordinates:
\begin{equation}\label{1.1}
  \left\{\begin{array}{l}
           v_t-u_x=0, \qquad\qquad\qquad\qquad\qquad\qquad\qquad\quad (x,t)\in\mathbb{R}^+\times\mathbb{R}^+,\\
           u_t+p(v)_x=\mu\left(\displaystyle\frac{u_x}{v}\right)_x+\kappa\left(\displaystyle\frac{-v_{xx}}{v^5}+\frac{5v_x^2}{2v^6}\right)_x,\quad\, (x,t)\in\mathbb{R}^+\times\mathbb{R}^+,
         \end{array}\right.
\end{equation}
with the initial and boundary conditions:
\begin{equation}\label{1.2}
  \left\{\begin{array}{l}
    (v,u)(x,0)=(v_0,u_0)(x),\quad\, x\geq0, \\
    u(0,t)=0,\quad v_x(0,t)=0,\quad t\geq0, \\
    (v,u)(+\infty,t)=(v_+,u_+),\quad\, t\geq0.
  \end{array}\right.
\end{equation}
Here the unknown functions  are the specific volume $v(t,x)>0$,  the velocity $u(t,x)$ and the pressure  $p=p(v(t,x))$ of the fluids.
The constants $\mu>0$ and $\kappa>0$ are the viscosity coefficient and the capillary coefficient, respectively, and $v_{+}>0$, $u_{+}$ are given constants.

Throughout this paper, we assume that the pressure $p(v)$ is a positive smooth function for $v>0$  satisfying
\begin{eqnarray}\label{1.3}
p^\prime(v)<0,\quad p^{\prime\prime}(v)>0, \quad \forall \,\,v>0,
\end{eqnarray}
and the initial data $u_0(x)$ satisfies
\begin{equation}\label{1.3-1}
u_0(0)=0
\end{equation}
as a compatibility condition.

System $(\ref{1.1})$ can be used to describe the motion of the compressible isothermal viscous fluids with internal capillarity,
see the pioneering  works by Van der Waals\cite{J. D. Van der Waals}, Korteweg \cite{Korteweg-1901} and  Dunn and Serrin \cite{J. E. Dunn-J. Serrin-1985}
for the derivations on the compressible Navier-Stokes-Korteweg system.   Notice that when $\kappa=0$,
the system $(\ref{1.1})$ is reduced to the compressible Navier-Stokes system.

There have been extensive studies on the mathematical aspects of the compressible Navier-Stokes-Korteweg system. For the initial value problem,
Hattori and Li \cite{H. Hattori-D. Li-1994, H. Hattori-D. Li-1996, H. Hattori-D. Li-1996-2} proved  the local existence and  global existence of
smooth solutions around constant states
in the Sobolev space  $H^s(\mathbb{R}^n)$ with $s\geq4+\frac{n}{2}$ and $n=2,3$. Since then, many authors study the global existence and  large-time behavior of strong solutions in the Sobolev space
framework, cf. \cite{Chen-2012-1, Chen-2012-2, Chen-2012-3, Chen-2012-4,Chen-2012-5,F-Charve-B-Haspot-2011, Y.-P. Li-2012, Wang-2015, Y. J. Wang-Z. Tan-2011}
and the references therein.
Danchin and Desjardins \cite{R. Danchin-B. Desjardins-2001} and Haspot\cite{B. Haspot-2009, B. Haspot-2016} obtained the global existence
and uniqueness of strong  solutions in  Besov space. Concerning the global existence of  weak solutions to the initial  value problem of the Korteweg system,
the readers are  referred   to \cite{R. Danchin-B. Desjardins-2001, B. Haspot-2011}  for small initial data  in the whole space $\mathbb{R}^2$,
\cite{D. Bresch-B. Desjardins-C. K. Lin-2003} for  large initial data  in a periodic domain $\mathbb{T}^d$ with $d=2$ or $3$, and \cite{B. Haspot-2011,Germain-LeFloch-2012} for  large initial data in the whole space $\mathbb{R}$.

 For the initial-boundary value problem,  fewer results have been obtained so far.
 Kotschote \cite{M. Kotschote-2008,M. Kotschote-2010} proved the local well-posedness  of strong solutions to the compressible Korteweg system in a bounded domain of $\mathbb{R}^n$ with $n\geq1$.
 Later, the global  existence  and  exponential decay of strong solutions to the Korteweg system  in a bounded domain of $\mathbb{R}^n$ ($n\geq1$) were obtained by Kotschote in \cite{E-Tsyganov-2008}.
 Tsyganov \cite{M. Kotschote-2014} discussed  the global existence and time-asymptotic behavior of  weak solutions for an isothermal system on the interval $[0,1]$. The time-asymptotic profiles in both \cite{M. Kotschote-2014} and \cite{E-Tsyganov-2008} are non-constant stationary solutions.

In the one-dimensional case,   the time-asymptotic nonlinear stability of some elementary waves
(such as rarefaction wave,  viscous shock wave and viscous contact discontinuity, etc.) to the initial value problem of the compressible Navier-Stokes-Korteweg system has been studied in  \cite{Chen-2012-1, Chen-2012-2, Chen-2012-3, Chen-2012-5}
However, to the best of our knowledge,  there is no result on the larger-time behavior of solutions for the initial-boundary value problem of the compressible Korteweg system toward these wave-like profiles.
The main purpose of this paper is devoted to this problem and as a first step  to this goal,
we will show the large-time behavior of strong solutions of  problem  (\ref{1.1})-(\ref{1.2}) toward the viscous shock profiles.
More precisely, we shall prove that the asymptotic profile of  (\ref{1.1})-(\ref{1.2}) is a viscous shock profile, which is suitably away from the boundary, and if the initial data is a small  perturbation  of the shifted viscous shock profile and the strength of the shifted viscous shock profile is sufficiently small,
then the problem (\ref{1.1})-(\ref{1.2}) has   a unique global (in time) strong solution, which tends  toward the shifted viscous shock profile as time goes to infinity.
The precise statement of our main result can be found in Theorem 2.1 below.

The asymptotic nonlinear stability of viscous shock profiles for various equations from  fluid mechanics in the half space $\mathbb{R}^+$ with boundary effect has been investigated by many authors. Liu and Yu \cite{T.-P.-Liu-S.-H.-Yu-1997} first studied the stability of the viscous shock profile for the Burgers equation with a Dirichlet boundary condition. Then  such a result was extended by Liu and Nishihara \cite{T-P liu-K. Nishihara-1997}  to the generalized Burgers equation. For the case of system, Matsumura and Mei \cite{Matsumura-Mei}  considered the  stability of viscous shock wave to  the p-system with viscosity and a Dirichlet boundary condition. Huang, Matsumura and Shi \cite{Huang-2003} obtained the nonlinear stability of  viscous shock  wave for an inflow problem of the isentropic compressible Navier-Stokes equations. The stability of a supposition of two viscous shock waves  for an initial-boundary value problem of the full compressible Navier-Stokes system   was shown by  Huang and Matsumura \cite{Huang-2009}.  For the other  results on the asymptotic stability of solutions with boundary effect,  we refer to \cite{T-P liu-A. Matsumura-K. Nishihara-1998,A. Matsumura-K. Nishihara-1986,P-LIU-Nishihara-1999,K. Nishihara-T.Yong-1999,Huang-2012,Huang-2009,Huang-2003} and the references therein. The present paper is motivated by the works \cite{Matsumura-Mei} and \cite{Chen-2012-3}. We first use the ideas of \cite{Matsumura-Mei} to  identify the asymptotic profile for the original problem (\ref{1.1})-(\ref{1.2}) and then show the time-asymptotic stability of this profile by  energy method similar to \cite{Chen-2012-3}. Compared with the case of \cite{Matsumura-Mei}, there are two additional  difficulties  in our analysis.  First,  some new difficulty boundary terms, such as $|\int_0^t(\frac{\phi_{xx}\psi_{xx}}{v^5})|_{x=0}d\tau|$ and $|\int_0^t(\frac{\phi_{xxx}\phi_{txx}}{v^5})|_{x=0}d\tau|$,  appear  in the process of energy estimates  due to  the Korteweg tensor  and the   estimates of the unknown functions $\psi_{xx}$ and $\phi_{xxx}$. We control these boundary terms by making use of the estimates of $A(t)$ (see (\ref{3.19})), the first equation of system (\ref{3.3}) and the  a priori assumption  (\ref{3.13}) repeatedly (see the proof of (\ref{3.17})-(\ref{3.18-3}) for details). The second difficult is to the estimate the nonlinear term $F$ (see (\ref{3.4})), which is more complicated than that of  \cite{Matsumura-Mei}. We overcome this difficulty by some delicate analysis. The main result of this paper can be viewed as an extension of \cite{Matsumura-Mei} to the case of  compressible fluid with  capillarity effect.

The layout of this paper is   as follows.  In Sections 2,  we first recall  the existence and properties of viscous shock profiles to system  (\ref{1.1}) in the whole space $\mathbb{R}$.
Then we identify the appropriate asymptotic profile for the  original problem (\ref{1.1})-(\ref{1.2})
and state  our main result.  In Section 3, we reformulate the original problem,  and then  perform the a priori estimates to the reformulated system.  The proof of our main result is given at the end of this section.

{\bf Notations:} Throughout this paper, $C$  denotes some generic constant which may vary in
different estimates. If the dependence needs to be explicitly pointed
out, the notation $C(\cdot,\cdots,\cdot)$ or $C_i(\cdot,\cdots,\cdot)(i\in {\mathbb{N}})$ is used.
 For function spaces, $L^p(\mathbb{R}^+)$ ($1\leq p\leq+\infty$) denotes   the standard  Lebesgue
space with the norm
 \[\|f\|_{L^p(\mathbb{R}^+)}=\left(\int_0^\infty\left|f(x)\right|^pdx\right)^{\frac{1}{p}}.\]
 $H^l(\mathbb{R}^+)$ is the
usual $l$-th order Sobolev space with its norm
$$\|f\|_{l}=\left(\sum_{i=0}^{l}\|\partial_x^if\|^2\right)^{\frac{1}{2}}\quad with\quad \|\cdot\|\triangleq\|\cdot\|_{L^2(\mathbb{R}^+)}.$$

\section{Preliminaries and main result}
\setcounter{equation}{0}
In this section,  we shall present  our main result on the time-asymptotic behavior of solutions to the initial-boundary value problem (\ref{1.1})-(\ref{1.2}).
To  this end, we first recall the existence and properties of viscous shock profiles to system  (\ref{1.1}) in the whole space $\mathbb{R}$,
 and then identify the appropriate asymptotic profiles for the  original problem. Finally, we give the main result of this paper  at the end of this section.
\subsection{ Viscous shock profiles}

Viscous shock profiles are the  traveling wave solutions of system  (\ref{1.1}) in the whole space $\mathbb{R}$ of the form:
\begin{equation}\label{2.1}
  (v,u)(x,t)=(V,U)(\xi),\quad \xi=x-st,
\end{equation}
which satisfies
\begin{eqnarray}\label{2.2}
\left\{\begin{array}{ll}
          -sV_\xi-U_\xi=0,\\[2mm]
          \displaystyle-sU_\xi+p(V)_\xi=\mu\left(\frac{U_\xi}{V}\right)_\xi+\kappa\left(\frac{-V_{\xi\xi}}{V^5}+\frac{5V_\xi^2}{2V^6}\right)_\xi
 \end{array}\right.
\end{eqnarray}
with the boundary conditions
\begin{equation}\label{2.3}
 (V, U)(\xi)\longrightarrow(v_{\pm}, u_{\pm}),\quad as\,\,\xi \rightarrow\pm\infty,
\end{equation}
where $s$ is the wave speed,  $\xi$ is the traveling wave variable and $v_{\pm}>0, u_{\pm}$ are given constants.

Integrating (\ref{2.2}) with respect to $\xi$ once, we obtain
\begin{eqnarray}\label{2.4}
\left\{\begin{array}{ll}
          sV+U=a_1,\\[2mm]
       \displaystyle\mu\frac{U_\xi}{V}-\kappa\frac{V_{\xi\xi}}{V^5}+\frac{5\kappa V_\xi^2}{2V^6}=-sU+p(V)+a_2,
 \end{array}\right.
\end{eqnarray}
where $a_1,a_2$ are constants satisfying
\begin{eqnarray}\label{2.5}
\left\{\begin{array}{ll}
          a_1=sv_{+}+u_{+}=sv_{-}+u_{-},\\[2mm]
       a_2=su_{+}-p(v_{+})=su_{-}-p(v_{-}),
 \end{array}\right.
\end{eqnarray}
which shows that $s$ and $v_{\pm},u_{\pm}$ satisfy the Rankine-Hugoniot condition for system (1.1):
\begin{eqnarray}\label{2.6}
\left\{\begin{array}{ll}
          s(v_{+}-v_{-})=(u_{-}-u_{+}),\\[2mm]
       s(u_{+}-u_{-})=p(v_{+})-p(v_{-}).
 \end{array}\right.
\end{eqnarray}
 Moreover, (\ref{2.6}) can be reduced to
\begin{eqnarray}\label{2.7}
s^2=\frac{p(v_{-})-p(v_{+})}{v_+-v_-}.
\end{eqnarray}
In this paper, we only consider the case $s>0$, i.e., the 2-shock wave, and the results for the  case $s<0$ follow similarly.   The Lax entropy condition for the 2-shock  (cf. \cite{P.D.Lax-1957}) is
\begin{equation}\label{2.8}
\sqrt{-p^\prime(v_+)}=\lambda_2(v_+)<s<\lambda_2(v_-)=\sqrt{-p^\prime(v_-)},
\end{equation}
then it follows from (\ref{1.3}), $(\ref{2.5})_1$ and (\ref{2.7}) that
\begin{equation}\label{2.9}
0<v_-<v_+,\quad u_+<u_-.
\end{equation}

We assume that $u_-=0$ throughout  this paper. Then for any given $(v_+,u_+)$ with  $v_+>0$ and $u_+<0$, the constants $v_-(0<v_-<v_+)$ and $s>0$ can be uniquely determined by the Rankine-Hugoniot condition (\ref{2.6}).

The existence and properties of the traveling wave solutions to system (\ref{1.1}) are summarized in the following theorem.
\begin{Proposition}(\cite{Chen-2012-3})
Let (\ref{1.3}) and (\ref{2.9}) hold. If
\[
\frac{\mu^2s^2v_-^8}{\kappa}-\left(\frac{10v_+}{v_-}-6\right)v_+^5\left(p^\prime(v_+)+s^2\right)>0,
\]
then there exists a monotone viscous shock profile $(V, U)(x-st)$ to system (\ref{1.1}), which is unique up to a shift and satisfies $V_\xi>0$, $U_\xi<0$, and
\begin{equation}\label{2.11}
  \left\{\begin{array}{l}
           |V(\xi)-v_-|\le C\delta e^{-c_1|\xi|},\quad |U(\xi)-u_-|\le C\delta e^{-c_1|\xi|},\quad \forall \xi\le0, \\
           |V(\xi)-v_+|\le C\delta e^{-c_1|\xi|},\quad|U(\xi)-u_+|\le C\delta e^{-c_1|\xi|},\quad \forall \xi\ge0, \\
           \displaystyle \left|\frac{d^k}{d\xi^k}V(\xi)\right|+\left|\frac{d^k}{d\xi^k}U(\xi)\right|\le C\delta^2 e^{-c_1|\xi|},\quad \forall\,\xi\in\mathbb{R},\forall\, k\ge1,
         \end{array}\right.
\end{equation}
where $\delta:=v_+-v_-$, and $c_1, C$ are two positive constants depending only on $v_+,v_-, s, \mu$ and $\kappa$.
\end{Proposition}

\subsection{Asymptotic profile for the original problem}
For the initial-boundary value problem (\ref{1.1})-(\ref{1.2}),  there exists  an initial boundary layer $(u(x,t)-U(x-st))|_{(x,t)=(0,0)}=u_--U(0)=-U(0)$ since $U(0)$ is always less than $u_-=0$. Consequently, the solutions of the original problem (\ref{1.1})-(\ref{1.2}) may not converge to a shifted viscous shock profile $(V,U)(x-st+\alpha)$ with the shift $\alpha$ determined by the initial data and the viscous shock profile. In order to get convergence, as \cite{Matsumura-Mei}, we consider the solution of the original problem (\ref{1.1})-(\ref{1.2}) in a neighborhood of $(V,U)(x-st+\alpha-\beta)$, where $\alpha$ is a shift to be determined later, and $\beta\gg1$ is a constant such that the initial boundary layer around the shifted wave $\left|(u(x,t)-U(x-st-\beta)\right|_{(x,t)=(0,0)}=\left|(u_--U(-\beta)\right|=\left|U(-\beta)\right|\ll1$. Now we determine the shift $\alpha$ for some given constant $\beta\gg1$.

First, it follows from $(\ref{1.1})_1$ and $(\ref{2.2})_1$ that
\begin{equation} \label{2.12}
\left(v-V\right)_t=\left(u-U\right)_x, \quad (V,U)=(V,U)(x-st+\alpha-\beta).
\end{equation}
Integrating (\ref{2.12}) with respect to $t$ and $x$ over $[0,t]\times\mathbb{R}^+$ and using the boundary condition $(\ref{1.2})_2$, we have
\begin{equation} \label{2.13}
\displaystyle\int_0^\infty[v(x,t)-V(x-st+\alpha-\beta)]dx
=\displaystyle\int_0^\infty[v_0(x)-V(x+\alpha-\beta)]dx+\int_0^tU(-s\tau+\alpha-\beta)d\tau.
\end{equation}

By the idea of  conservation of mass principle, we suppose that
\begin{equation} \label{2.14}
\displaystyle\int_0^\infty[v(x,t)-V(x-st+\alpha-\beta)]dx\rightarrow0,\quad as\,\, t\rightarrow\infty.
\end{equation}
Then we set
\begin{equation} \label{2.15}
I(\alpha)=\displaystyle\int_0^\infty[v_0(x)-V(x+\alpha-\beta)]dx+\int_0^\infty U(-s\tau+\alpha-\beta)d\tau.
\end{equation}
We see from (\ref{2.13}) and (\ref{2.14}) that the shift $\alpha$ must satisfies $I(\alpha)=0$. Differentiating  (\ref{2.15}) with respect to $\alpha$, we have
\begin{equation} \label{2.16}
\aligned
I^\prime(\alpha)&=-\displaystyle\int_0^\infty V^\prime(x+\alpha-\beta)]dx+\int_0^\infty U^\prime(-s\tau+\alpha-\beta)d\tau\\
&=-v_++V(\alpha-\beta)+\frac{1}{s}U(\alpha-\beta)\\
&=-v_++\frac{1}{s}(sv_-+u_-)=v_--v_+,
\endaligned
\end{equation}
where we have used $(\ref{2.4})_1$ and the assumption that $u_-=0$. Integrating (\ref{2.16}) with respect to $\alpha$ over $[0,\alpha]$ yields
\begin{equation} \label{2.17}
\aligned
I(\alpha)&=I(0)+(v_--v_+)\alpha\\
&=\displaystyle\int_0^\infty[v_0(x)-V(x-\beta)]dx+\int_0^\infty U(-st-\beta)dt+(v_--v_+)\alpha.
\endaligned
\end{equation}
Since $I(\alpha)=0$, the shift $\alpha=\alpha(\beta)$ is determined explicitly by
\begin{equation} \label{2.18}
\alpha=\frac{1}{v_+-v_-}\left\{\displaystyle\int_0^\infty[v_0(x)-V(x-\beta)]dx+\int_0^\infty U(-st-\beta)dt\right\}.
\end{equation}
Then we deduce from (\ref{2.14}), (\ref{2.15}) and $I(\alpha)=0$ that
\begin{equation} \label{2.19}
\aligned
\displaystyle\int_0^\infty[v(x,t)-V(x-st+\alpha-\beta)]dx
&=I(\alpha)-\int_t^\infty U(-s\tau+\alpha-\beta)d\tau\\
&=-\int_t^\infty U(-s\tau+\alpha-\beta)d\tau\rightarrow0, \quad as\,\, t\rightarrow\infty,
\endaligned
\end{equation}
which, together with (\ref{2.13}) implies that
\begin{equation} \label{2.20}
\displaystyle\int_0^\infty[v_0(x)-V(x+\alpha-\beta)]dx=-\int_0^\infty U(-s\tau+\alpha-\beta)d\tau.
\end{equation}
By a similar argument to the second equation of (\ref{1.1}), we have
\begin{equation}\label{2.21}
  \begin{split}
     & \int_0^{\infty}\left[u_0(x)-U(x+\alpha-\beta)\right]dx+\int_0^{\infty}\left[p(v(0,t))-p(V(-st+\alpha-\beta))\right] dt \\
     - & \mu\int_0^{\infty}\left[\frac{u_x(0,t)}{v(0,t)}-\frac{U^{\prime}(-st+\alpha-\beta)}{V(-st+\alpha-\beta)}\right] dt\\
     +& \kappa\int_0^{\infty}\left[\frac{v_{xx}(0,t)}{v(0,t)^5}-\frac{5v_x^2(0,t)}{2v^6(0,t)}-\frac{V^{\prime\prime}(-st+\alpha-\beta)}{V^5(-st+\alpha-\beta)}+\frac{5}{2}\frac{(V^{\prime}(-st+\alpha-\beta))^2}{V(-st+\alpha-\beta)}\right] dt=0.
  \end{split}
\end{equation}
Here we expect that $v(0,t)$, $u_x(0,t)=v_t(0,t)$ and $v_{xx}(0,t)$  can be controlled  by the effects of boundary, viscosity and capillarity such that   (\ref{2.21}) holds with the same shift $\alpha$ defined by (\ref{2.18}). This is possible because $v(0,t)$ is not specified.

Thus, by the  above heuristical analysis, we expect the asymptotic profile for the original problem  (\ref{1.1})-(\ref{1.2}) is the shifted traveling waves $(V,U)(x-st+\alpha-\beta)$ with the shift $\alpha=\alpha(\beta)\ll1$ and the constant $\beta\gg1$.

\subsection{Main result}

To state our main result,  we suppose that for some  constant $\beta>0$,
\begin{equation}\label{2.22}
v_0(x)-V(x-\beta)\in H^1(\mathbb{R}^+)\cap L^1(\mathbb{R}^+),\quad u_0(y)-U(y-\beta)\in H^1(\mathbb{R}^+)\cap L^1(\mathbb{R}^+).
\end{equation}
Set
\begin{equation}\label{2.23}
 \left(\Phi_0,\Psi_0\right)(x)=-\int_x^\infty\left(v_0(y)-V(y-\beta), u_0(y)-U(y-\beta)\right)dy,
\end{equation}
and assume that
\begin{equation}\label{2.24}
(\Phi_0,\Psi_0)\in L^2(\mathbb{R}^+).
\end{equation}
Then our  main result is stated as follows.
\begin{Theorem}
 Assume that  $v_+>0$, $u_+<0$ and $\beta>0$ are  constants.   Let (\ref{1.3-1}), (\ref{2.22})-(\ref{2.24}) and the conditions  listed in Proposition 2.1 hold,  and $(V,U)(x-st)$ be the viscous shock  profile obtained in Proposition 2.1. Then there exists two small positive constants $\delta_0$ and $\varepsilon_0$ such that if $0<\delta:=v_+-v_-\leq\delta_0$ and
 \begin{equation}\label{2.25}
 \|\Phi_0\|_3+\|\Psi_0\|_2+\beta^{-1}\leq\varepsilon_0,
\end{equation}
 then the initial-boundary value problem (\ref{1.1})-(\ref{1.2}) admits a unique global solution $(u,v)(t,x)$ satisfying
\begin{equation}\label{2.26}
 \aligned
v(x,t)-V(x-st+\alpha-\beta)\in C\left([0,\infty; H^2(\mathbb{R}^+)\right)\cap L^2\left([0,\infty; H^3(\mathbb{R}^+)\right),\\
u(x,t)-U(x-st+\alpha-\beta)\in C\left([0,\infty; H^1(\mathbb{R}^+)\right)\cap L^2\left([0,\infty; H^2(\mathbb{R}^+)\right),
\endaligned
\end{equation}
where $\alpha=\alpha(\beta)$ is the  shift determined by (\ref{2.18}). Moreover, the following asymptotic behavior of solutions holds:
\begin{equation}\label{2.27}
 \displaystyle \lim_{t\rightarrow\infty}\sup_{x\in \mathbb{R}^+}\left|(v,u)(x,t)-(V,U)(x-st+\alpha-\beta)\right|=0.
\end{equation}
\end{Theorem}

\section{Proof of Theorem 2.1}
\setcounter{equation}{0}
This section is devoted to proving Theorem 2.1. We begin with a reformulation of the original problem (\ref{1.1})-(\ref{1.2}).
\subsection{Reformulation of the original problem}
First, we define the perturbation functions $(\phi,\psi)(t,x)$ as follows:
\begin{equation}\label{3.1}
(\phi, \psi)(x,t)=\displaystyle-\int_x^{\infty}\left(v(y,t)-V(y-st+\alpha-\beta), u(y,t)-U(y-st+\alpha-\beta)\right) dy
\end{equation}
for $(x,t)\in \mathbb{R}^+\times  \mathbb{R}^+$. Then, we have
\begin{equation}\label{3.2}
 (v,u)(x,t)=(\phi_x(x,t)+V(x-st+\alpha-\beta),\,\,\psi_x(x,t)+U(x-st+\alpha-\beta)).
\end{equation}
Substituting (\ref{3.2}) into (\ref{1.1}), using (\ref{2.2}) and integrating the system in $x$ over $[x,+\infty)$,  we obtain
\begin{equation}\label{3.3}
 \left\{\begin{array}{l}
          \phi_t-\psi_x=0, \\
          \psi_t+\displaystyle p^{\prime}(V)\phi_x-\frac{\mu}{V}\psi_{xx}+\frac{\kappa\phi_{xxx}}{v^5}+\mu\frac{U^{\prime}\phi_x}{V^2}=F,
        \end{array}\right.
\end{equation}
where
\begin{equation}\label{3.4}
\aligned
    F=&-\{p(v)-p(V)-p^{\prime}(V)\phi_x\}-\frac{\mu\psi_{xx}\phi_x}{vV}+\kappa V^{\prime\prime}\left(\frac{1}{V^5}-\frac{1}{v^5}\right)\\
      &+\frac{5\kappa}{2}\frac{\phi_{xx}^2+2\phi_{xx}V^{\prime}}{v^6}+\frac{5\kappa}{2}(V^{\prime})^2\left(\frac{1}{V^6}-\frac{1}{v^6}\right)+\frac{\mu U^{\prime}\phi_x^2}{vV^2}.
  \endaligned
\end{equation}
The initial data of system (\ref{3.3}) satisfy
\begin{equation}\label{3.5}
 \aligned
\phi_0(x)&:=\phi(x,0)
      =-\int_x^{\infty}[v_0(y)-V(y-\beta)]\ dy+\int_x^{\infty}[V(y+\alpha-\beta)-V(y-\beta)]\ dy\\
      &=\Phi_0(x)+\int_x^{\infty}[V(y+\alpha-\beta)-V(y-\beta)]\ dy\\
      &=\Phi_0(x)+\int_x^{\infty}\int_0^{\alpha}V^{\prime}(y+\theta-\beta)\ d\theta dy\\
      &=\Phi_0(x)+\int_0^{\alpha}\int_x^{\infty}V^{\prime}(y+\theta-\beta)\ dy d\theta \\
      &=\Phi_0(x)+\int_0^{\alpha}[v_+-V(x+\theta-\beta)]\ d\theta,
 \endaligned
\end{equation}
and
\begin{equation}\label{3.6}
 \aligned
     \psi_0(x)&:=\psi(x,0) =-\int_x^{\infty}[u_0(y)-U(y+\alpha-\beta)]\ dy \\
      & =-\int_x^{\infty}[u_0(y)-U(y-\beta)]\ dy+\int_x^{\infty}[U(y+\alpha-\beta)-U(y-\beta)]\ dy\\
      &=\Psi_0(x)+\int_0^{\alpha}[u_+-U(x+\theta-\beta)]\ d\theta,
  \endaligned
\end{equation}
where the functions $(\Phi_0,\Psi_0)(x)$ are defined in (\ref{2.23}).

Using (\ref{3.1}), (\ref{3.2}) and (\ref{2.19}), we have
\begin{equation}\label{3.7}
  \begin{split}
    \phi(0,t)&=-\int_0^{\infty}[v(y,t)-V(y-st+\alpha-\beta)]\ dy \\
           &=\int_t^{\infty}U(-s\tau+\alpha-\beta)\ d\tau:=A(t),
  \end{split}
\end{equation}
\begin{equation}\label{3.8}
  \begin{split}
    \psi_x(0,t)&=u(0,t)-U(-st+\alpha-\beta) \\
      &=u_--U(-st+\alpha-\beta)\\
      &=-U(-st+\alpha-\beta)=A^{\prime}(t).
  \end{split}
\end{equation}
On the other hand,
 it follows from $v_x(0,t)=0$ that
\[
  \aligned
    0&=v_x|_{x=0}=\phi_{xx}|_{x=0}+V^{\prime}(x-st+\alpha-\beta)|_{x=0} \\
      &=\phi_{xx}(0,t)+V^{\prime}(-st+\alpha-\beta),
  \endaligned
\]
consequently,
\begin{equation}\label{3.9}
    \aligned
     \phi_{xx}(0,t)&=-V^{\prime}(-st+\alpha-\beta)=-\frac{1}{s}(sV^{\prime})(-st+\alpha-\beta) \\
      &=\frac{1}{s}U^{\prime}(-st+\alpha-\beta)=\frac{1}{s^2}A^{\prime\prime}(t),
  \endaligned
\end{equation}
where we have used  $(\ref{2.2})_1$ and (\ref{3.8}). We assume that  $\phi_0(0)=\phi(0,0)=A(0)$ as a compatible  condition.

In what follows, we seek the solutions of the initial-boundary  value problem (\ref{3.4})-(\ref{3.9}) in the following function space:
\begin{equation}\label{3.10}
X_{M}(0, T)=\left\{(\phi, \psi)(x,t)\left|
\begin{array}{c}
\phi(x,t)\in C(0, T; H^3(\mathbb{R})),\phi_x(x,t)\in L^2(0, T; H^{3}(\mathbb{R})),\\[2mm]
\psi(x,t)\in C(0, T; H^{2}(\mathbb{R})),\psi_x(x,t)\in L^2(0, T; H^{2}(\mathbb{R})),\\[2mm]
\displaystyle\sup_{t\in[0, T]}\left\{\|\phi(t)\|_3+\|\psi(t)\|_2\right\}\leq M,
\end{array}
\right.\right\}
\end{equation}
where $0\leq T\leq\infty$, and $M$ is some positive constant.

For the problem (\ref{3.4})-(\ref{3.9}), we have the following theorem, which leads to  Theorem 2.1 immediately.
\begin{Theorem}
 Suppose  that the conditions of Theorem 2.1 hold, then there exists positive constants $\delta_0$, $\varepsilon_1$ and $C_0$ which are independent of the time $t$ and the initial data $\phi_0,\psi_0$ such that if $0<\delta:=v_+-v_-\leq\delta_0$ and $
 \|\phi_0\|_3+\|\psi_0\|_2+\beta^{-1}\leq\varepsilon_1$, then the initial-boundary  value problem (\ref{3.4})-(\ref{3.9}) has a unique global solution $(\phi,\psi)(t,x)\in X_{\hat{M}}[0,\infty)$ with $\hat{M}=2\sqrt{C_0(\|\phi_0\|_3^2+\|\psi_0\|_2^2+e^{-c_1\beta})}$.
 Moreover, it holds that
 \begin{equation}\label{3.11}
\|\phi(t)\|_{3}^2+\|\psi(t)\|_{2}^2+\int_0^t\|(\phi_x,\psi_x)(\tau)\|_{3,2}^2\ d\tau
     \le  C_0(\|\phi_0\|_3^2+\|\psi_0\|_2^2+e^{-c_1\beta})
\end{equation}
for any $t\in[0,+\infty)$, and
\begin{equation}\label{3.12}
 \displaystyle \lim_{t\rightarrow\infty}\sup_{x\in \mathbb{R}^+}\left|(\phi_x, \psi_x)(x,t)\right|=0.
\end{equation}
\end{Theorem}

Notice that the initial conditions in Theorems 3.1 and 2.1 are different. The following lemma state the relation between them.
\begin{Lemma}
Suppose that the condition (\ref{2.22}) holds. Then we have
\begin{itemize}
\item [(i)] The shift $\alpha\rightarrow0$ if $\|\Phi_0\|_2\rightarrow0$ and $\beta\rightarrow\infty$;
\item [(ii)] $\|\phi_0\|_3+\|\psi_0\|_2\rightarrow0$ if $\|\Phi_0\|_3+\|\Psi_0\|_2\rightarrow0$ and $\beta\rightarrow\infty$.
\end{itemize}
\end{Lemma}

Since the proof of Lemma 3.1 (i) and (ii) are  almost the same as those  of Lemma 2.1 and Lemma 3.2 in \cite{Matsumura-Mei} respectively, we omit the details here for  brevity.

Theorem 3.1 will be obtained by combining  the following local existence and a priori estimates.
\begin{Proposition} [Local existence]
 Under the assumptions of Theorem 3.1, suppose that $\|\phi_0\|_3+\|\psi_0\|_2\leq M$, then there exists a positive constant $t_0$ depending only on $M$ such that the initial-boundary value  problem (\ref{3.4})-(\ref{3.9}) admits a unique  solution $(\phi, \psi)(x,t)\in X_{2M}(0,t_0)$.
\end{Proposition}

\begin{Proposition} [A priori estimates]
Suppose that $(\phi,\psi)(x,t)\in X_{M}(0, T)$ is a solution of the initial-boundary value  problem (\ref{3.4})-(\ref{3.9}) obtained in Proposition 3.1 for some  positive constants $T$ and $M$. Then there exists three positive constants $\delta_0\ll1$, $\varepsilon_1\ll1$ and $C_0$ which are independent of $T$ such that if $0<\delta\leq\delta_0$ and
\begin{equation}\label{3.13}
 N(T):=\displaystyle\sup_{t\in[0, T]}\left\{\|\phi(t)\|_3+\|\psi(t)\|_2\right\}\leq \varepsilon
\end{equation}
 with the constant $\varepsilon$ satisfying $0<\varepsilon\leq\varepsilon_1$, the solution $(\phi,\psi)$ of the problem (\ref{3.4})-(\ref{3.9}) satisfies the estimates (\ref{3.11}) for all $t\in[0,T]$.
\end{Proposition}

 Proposition 3.1 can be proved  by using the dual argument and iteration technique, which is similar to that of Theorem 1.1 in \cite{H. Hattori-D. Li-1994}, the details are omitted here. The proof of Proposition 3.2 is more subtle and will be given in the next subsection.

\subsection{A priori estimates}
Due to Lemma 3.1 (i) and the assumptions in Theorem 3.1 and Proposition 3.2, without loss of generality,  we assume that $|\alpha|<1,\beta>1$ and $N(t)<1$ hereafter.

Let $(\phi,\psi)(x,t)\in X_{M}(0, T)$ be a solution of the initial-boundary value  problem (\ref{3.4})-(\ref{3.9})  for some  positive constants $T$ and $M$. Then we have from Proposition 2.1, the Sobolev inequality $\|f\|_{L^\infty}\leq \|f\|_1$ for $f(x)\in H^1$, and the smallness of  $\varepsilon_1>0$ that
\begin{equation}\label{3.14}
   v(x,t)=V(x,t)+\phi_x(x,t)\le v_++\|\phi_x\|_{L^{\infty}}\le v_++\|\phi_x(t)\|_1\le \frac{3}{2}v_+, \quad \forall \,(x,t)\in [0,T]\times\mathbb R^+,
\end{equation}
\begin{equation}\label{3.15}
 v(x,t)=V(x,t)+\phi_x(x,t)\ge v_--\|\phi_x\|_{L^{\infty}}\ge v_--\|\phi_x(t)\|_1\ge \frac{v_-}{2},\quad \forall \,(x,t)\in [0,T]\times\mathbb R^+.
\end{equation}

We first give the following boundary estimates.

\begin{Lemma}[Boundary estimates] Under the assumptions of Proposition 3.2, there exists a positive constant $C>0$ which is independent of $\alpha, \beta$ and $T$  such that the following boundary estimates:
\begin{equation}
\label{3.16}
\aligned
&\left|\int_0^t(\phi\psi)|_{x=0}d\tau\right|+\left|\int_0^t(\psi\psi_x)|_{x=0}d\tau\right|+\left|\int_0^t(\psi_{x}\phi_x)|_{x=0}d\tau\right|+\left|\int_0^t(\psi_t\psi_x)|_{x=0}d\tau\right|\\
&+\left|\int_0^t(\psi_{x}\psi_{xx})|_{x=0}d\tau\right|
+\left|\int_0^t(\psi_{xt}\psi_{xx})|_{x=0}d\tau\right|\leq Ce^{-c_1\beta},
\endaligned
\end{equation}
\begin{eqnarray}
\label{3.17}
&&\left|\int_0^t\left(\frac{\phi_{xx}\psi}{p^\prime(V)v^5}\right)\bigg|_{x=0}d\tau\right|+\left|\int_0^t\left(\frac{\phi_x\phi_{xx}}{v^5}\right)\bigg|_{x=0}d\tau\right|
+\left|\int_0^t\left(\frac{\phi_{xx}\psi_{xx}}{v^5}\right)\bigg|_{x=0}d\tau\right|\nonumber\\
&&+\left|\int_0^t(\psi_{t}\phi_{xx})|_{x=0}d\tau\right|\leq Ce^{-c_1\beta},
\end{eqnarray}
\begin{equation}\label{3.18}
    \left|\int_0^t\left(\frac{\phi_{xxx}\phi_{txx}}{v^5}\right)\bigg|_{x=0}d\tau\right|
     \le C\varepsilon^{\frac{1}{2}}\int_0^t\|\phi_{xxxx}(\tau)\|^2\ d\tau+C\varepsilon^{\frac{1}{2}}e^{-\frac{4}{3}c_1\beta},
\end{equation}
\begin{equation}\label{3.18-1}
    \left|\int_0^t\left(\psi_{tx}\phi_{xxx}\right)|_{x=0} d\tau\right|
     \le \eta\int_0^t\|\phi_{xxxx}(\tau)\|^2\ d\tau+C_{\eta}\int_0^t\|\phi_{xxx}(\tau)\|^2\ d\tau+C_{\eta}e^{-c_1\beta}.
\end{equation}

\begin{equation}\label{3.18-2}
   \left|\int_0^t\left(p^{\prime}(V)\phi_{xx}\phi_{xxx}\right)|_{x=0}\ d\tau\right|\le  \eta\int_0^t\|\phi_{xxxx}(\tau)\|^2\ d\tau+C_{\eta}\int_0^t\|\phi_{xx}(\tau)\|_1^2d\tau,
   \end{equation}

\begin{equation}
\label{3.18-3}
\left|\int_0^t \left(\psi_{xx}\psi_{xxx}\right)|_{x=0}d\tau\right|
      \le C\int_0^t\|\psi_{xx}(\tau)\|_1^2\ d\tau+Ce^{-c_1\beta}
\end{equation}
hold for all  $0\leq t\leq T$ , where $c_1>0$ is constant defined in Proposition 2.1, $\eta$ is a small positive constant and $C_\eta$ is a positive constant depending on $\eta$.
\end{Lemma}
\noindent{\bf Proof.}~~ First, we show that the boundary data (\ref{3.7}): $\phi(0,t)=A(t):=\int_t^{\infty}U(-s\tau+\alpha-\beta)\ d\tau$ satisfies that
\begin{equation}
\label{3.19}
\aligned
\left|\frac{d^k}{dt^k}A(t)\right|&\leq Ce^{-c_1\beta}e^{-c_1st},\quad k=0,1,2,3, \\
\|A\|_{W^{3,1}}&\leq Ce^{-c_1\beta}.
  \endaligned
\end{equation}
In fact, since $|-s\tau+\alpha-\beta|=s\tau+\beta-\alpha$ because $s>0$ and $\beta>\alpha$, it follows from Proposition 2.1 that
\[|U(-s\tau+\alpha-\beta)|\leq Ce^{-c_1|-s\tau+\alpha-\beta|}\leq Ce^{-c_1(\beta-\alpha)}e^{-c_1s\tau}\leq Ce^{-c_1\beta}e^{-c_1s\tau}.\]
 Consequently, we have
\begin{equation}
\label{3.20}
|A(t)|\leq Ce^{-c_1\beta}e^{-c_1st}.
\end{equation}
Similarly,  by differentiating $A(t)$ with respect to $t$ and using Proposition 2.1, we can also obtain that $A^\prime(t), A^{\prime\prime}(t)$ and $A^{\prime\prime\prime}(t)$ are all bounded by $Ce^{-c_1\beta}e^{-c_1st}$, thus the inequalities in (\ref{3.19}) hold.

On the other hand, by the Sobolev inequality,  we have
\begin{equation}
\label{3.21}
\aligned
\left|\psi(0,t)\right|&\leq \sup_{x\in R^+}|\psi(x,t)|\leq CN(T), \\
\left|\phi_x(0,t)\right|+\left|\psi_x(0,t)\right|&\leq \sup_{x\in R^+}\left(|\phi(x,t)|+|\psi(x,t)|\right)\leq CN(T).
  \endaligned
\end{equation}

 Now we give the proofs of (\ref{3.16})-(\ref{3.18-3}). For (\ref{3.16}), since its proof is  almost the same as that in Lemma 4.1 of \cite{Matsumura-Mei},
 we omit  the details here. Next, we focus  on deducing (\ref{3.17})-(\ref{3.18-3}),
 which are new boundary terms compared with the the case of compressible Navier-Stokes system (see Lemma 4.1 in \cite{Matsumura-Mei}).
 Using the boundary conditions (\ref{3.9}), the estimates (\ref{3.19}) and (\ref{3.21}),  we have
 \begin{equation}
\label{3.22}
  \aligned
\left|\int_0^t\left(\frac{\phi_{xx}\psi}{p^\prime(V)v^5}\right)\bigg|_{x=0}d\tau\right|&\leq C\int_0^t|\phi_{xx}\psi|_{x=0}\ d\tau\le C\int_0^t|A^{\prime\prime}(\tau)||\psi(0,\tau)|d\tau \\
      & \leq CN(T)\int_0^t|A^{\prime\prime}(\tau)|\ d\tau\le Ce^{-c_1\beta},
  \endaligned
\end{equation}
\begin{equation}
 \label{3.23}
 \aligned
   \left|\int_0^t\left(\frac{\phi_x\phi_{xx}}{v^5}\right)\bigg|_{x=0}\ d\tau\right|&\leq C\int_0^t|\phi_x(0,\tau)\phi_{xx}(0,\tau)|\ d\tau \\
      &\le C\int_0^t|\phi_{x}(0,\tau)||A^{\prime\prime}(\tau)|\ d\tau\\
      &\le CN(T)\int_0^t|A^{\prime\prime}(\tau)|\ d\tau\le Ce^{-c_1\beta}.
  \endaligned
\end{equation}
To estimate the other boundary terms in (\ref{3.17}),
 we have by making use of the relation $\phi_{tx}=\psi_{xx}$ (see $(\ref{3.3})_1$), integration by parts , Proposition 2.1 and (\ref{3.19})-(\ref{3.21}) that
\begin{eqnarray} \label{3.24}
    \left|\int_0^t\left(\frac{\phi_{xx}\psi_{xx}}{v^5}\right)\bigg|_{x=0} d\tau\right|&&=\left|\int_0^t\left(\frac{\phi_{xx}\phi_{tx}}{v^5}\right)\bigg|_{x=0}d\tau\right|
      =\left|\int_0^t\left(\left(\phi_{xx}\frac{\phi_{x}}{v^5}\right)_t-\left(\frac{\phi_{xx}}{v^5}\right)_t\phi_x\right)\bigg|_{x=0} d\tau\right|\nonumber\\
      &&=\left(\phi_{xx}\frac{\phi_{x}}{v^5}\right)(0,t)-\left(\frac{\phi_{xx}}{v^5}\phi_x\right)(0,0)-\int_0^t\left(\phi_{xxt}\frac{\phi_{x}}{v^5}+5\frac{\phi_{xx}v_t\phi_x}{v^6}\right)\bigg|_{x=0} d\tau\nonumber\\
     &&\le C(|\phi_{xx}(0,t)\phi_x(0,t)|+|\phi_{xx}(0,0)||\phi_x(0,0)|)\nonumber\\
      &&\quad+C\int_0^t\left|\phi_{xxt}(0,\tau)\phi_x(0,\tau)\right| d\tau+C\left|\int_0^t\left(\frac{\phi_{xx}v_t}{v^6}\phi_x\right)(0,\tau)d\tau\right|\nonumber\\
    &&\le C(|A^{\prime\prime}(t)||\phi_x(0,t)|+|A^{\prime\prime}(0)||\phi_x(0,0)|)+C\int_0^t|A^{\prime\prime\prime}(\tau)||\phi_x(0,\tau)|d\tau\nonumber\\
      &&\quad+C\left|\int_0^t\left(\frac{\phi_{xx}v_t}{v^6}\phi_x\right)(0,\tau)d\tau\right|\nonumber\\
      &&\le CN(T)e^{-c_1\beta}+C\left|\int_0^t\left(\frac{\phi_{xx}v_t}{v^6}\phi_x\right)(0,\tau)d\tau\right|,
\end{eqnarray}
and
\begin{equation}
 \label{3.25}
\begin{split}
  &\left|\int_0^t\left(\frac{\phi_{xx}v_t}{v^6}\phi_x\right)(0,\tau)\ d\tau\right|\\
   \le &\left|\int_0^t\frac{\phi_{xx}(0,\tau)}{v^6(0,\tau)}\phi_x(0,\tau)(U^{\prime}+\psi_{xx}(0,\tau)) d\tau\right|\\
\le & C\delta^2\int_0^t|\phi_{xx}(0,\tau)\phi_x(0,\tau)|\ d\tau+\left|\int_0^t\frac{\phi_{xx}(0,\tau)\psi_{xx}(0,\tau)}{v^5(0,\tau)} d\tau\right|\cdot\left\|\frac{\phi_x}{v}(0,t)\right\|_{L^{\infty}_t}\\
\le & C\delta^2N(T)\int_0^t|A^{\prime\prime}(\tau)| d\tau+CN(T)\left|\int_0^t\left(\frac{\phi_{xx}\psi_{xx}}{v^5}\right)\bigg|_{x=0} d\tau\right|\\
\le& C\delta^2e^{-c_1\beta}+C\varepsilon\left|\int_0^t\left(\frac{\phi_{xx}\psi_{xx}}{v^5}\right)\bigg|_{x=0} d\tau\right|.
  \end{split}
\end{equation}
Combining (\ref{3.24})-(\ref{3.25}), and using the smallness of $\varepsilon$ yields
\begin{equation} \label{3.26}
  \left|\int_0^t\left(\frac{\phi_{xx}\psi_{xx}}{v^5}\right)\bigg|_{x=0}\ d\tau\right|\le Ce^{-c_1\beta}.
\end{equation}
Similarly, we have
\begin{eqnarray}\label{3.27}
    \left|\int_0^t\left(\psi_t\phi_{xx}\right)|_{x=0} d\tau\right|&&=\left|\int_0^t\psi_t(0,\tau)\frac{1}{s^2}A^{\prime\prime}(\tau)\ d\tau\right|\nonumber\\
     &&= \frac{1}{s^2}\left|\int_0^t[(\psi(0,\tau)A^{\prime\prime}(\tau))_{\tau}-\psi(0,\tau)A^{\prime\prime\prime}(\tau)]d\tau\right|\nonumber\\
     &&=\frac{1}{s^2}\left|\psi(0,t)A^{\prime\prime}(t)-\psi(0,0)A^{\prime\prime}(0)-\int_0^t\psi(0,\tau)A^{\prime\prime\prime}(\tau)d\tau\right|\nonumber\\
     &&\le  CN(T) e^{-c_1\beta}+CN(T)\int_0^t|A^{\prime\prime\prime}(\tau)| d\tau\nonumber\\
     &&\le C e^{-c_1\beta}.
\end{eqnarray}
Thus (\ref{3.17}) follows from (\ref{3.22}), (\ref{3.23}), (\ref{3.26}) and (\ref{3.27}) immediately.

 On the other hand, we derive from the equalities $\phi_{txx}(0,t)=A^{\prime\prime\prime}(t)$,  $\psi_{tx}(0,t)=A^{\prime\prime}(t)$, the Sobolev inequality, the Cauchy inequality  and (\ref{3.19}) that
\begin{equation}\label{3.30}
  \begin{split}
    \left|\int_0^t\left(\frac{\phi_{xxx}\phi_{txx}}{v^5}\right)\bigg|_{x=0}d\tau\right| \le &C\int_0^t|\phi_{xxx}(0,\tau)||A^{\prime\prime\prime}(\tau)|\ d\tau  \\
     \le&C\int_0^t\|\phi_{xxx}(\tau)\|^{\frac{1}{2}}\|\phi_{xxxx}(\tau)\|^{\frac{1}{2}} |A^{\prime\prime\prime}(\tau)|\ d\tau\\
     \le&C\varepsilon^{\frac{1}{2}}\int_0^t\|\phi_{xxxx}(\tau)\|^2\ d\tau+C\varepsilon^{\frac{1}{2}}e^{-\frac{4}{3}c_1\beta},
  \end{split}
\end{equation}
\begin{eqnarray}\label{3.31}
    \left|\int_0^t\left(\psi_{tx}\phi_{xxx}\right)|_{x=0} d\tau\right|&&\le  \int_0^t|\phi_{xxx}(0,\tau)||A^{\prime\prime}(\tau)|\ d\tau \nonumber\\
     &&\le  \int_0^t\|\phi_{xxx}(\tau)\|^{\frac{1}{2}}\|\phi_{xxxx}(\tau)\|^{\frac{1}{2}}|A^{\prime\prime}(\tau)|\ d\tau \nonumber\\
     &&\le \eta\int_0^t\|\phi_{xxxx}(\tau)\|^2\ d\tau+C_{\eta}\int_0^t(\|\phi_{xxx}(\tau)\|^2+|A^{\prime\prime}(\tau)|^2)\ d\tau\nonumber\\
     &&\le \eta\int_0^t\|\phi_{xxxx}(\tau)\|^2\ d\tau+C_{\eta}\int_0^t\|\phi_{xxx}(\tau)\|^2\ d\tau+C_{\eta}e^{-c_1\beta},
\end{eqnarray}
\begin{equation}\label{3.31-1}
 \begin{split}
    \left|\int_0^t\left(p^{\prime}(V)\phi_{xx}\phi_{xxx}\right)|_{x=0}\ d\tau\right|\le & C\int_0^t|\phi_{xx}(0,\tau)||\phi_{xxx}(0,\tau)| d\tau \\
    \le & C\int_0^t\|\phi_{xx}(\tau)\|^{\frac{1}{2}}\|\phi_{xxx}(\tau)\|\|\phi_{xxxx}(\tau)\|^{\frac{1}{2}}d\tau\\
    \le &\eta\int_0^t\|\phi_{xxxx}(\tau)\|^2\ d\tau+C_{\eta}\int_0^t\|\phi_{xx}(\tau)\|_1^2 d\tau,
 \end{split}
\end{equation}
\begin{eqnarray}\label{3.32}
   \left|\int_0^t \left(\psi_{xx}\psi_{xxx}\right)|_{x=0}d\tau\right|&&= \left|\int_0^t \left(\psi_{xx}\phi_{txx}\right)|_{x=0}d\tau\right| \le C\int_0^t|\psi_{xx}(0,\tau)||A^{\prime\prime\prime}(\tau)|\,d\tau\nonumber\\
      &&\le C\int_0^t\|\psi_{xx}(\tau)\|^{\frac{1}{2}}\|\psi_{xxx}(\tau)\|^{\frac{1}{2}}|A^{\prime\prime\prime}(\tau)|\ d\tau\nonumber\\
      &&\le C\int_0^t\|\psi_{xx}(\tau)\|_1^2\ d\tau+C\int_0^t|A^{\prime\prime\prime}(\tau)|\ d\tau\nonumber\\
      &&\le C\int_0^t\|\psi_{xx}(\tau)\|_1^2\ d\tau+Ce^{-c_1\beta}.
\end{eqnarray}
From (\ref{3.30})-(\ref{3.32}), we obtain  (\ref{3.18})-(\ref{3.18-3}) respectively. This completes the proof of Lemma 3.2.

Now we give the $L^2$-estimates on $(\phi,\psi)(x,t)$.

\begin{Lemma}
 Under the assumptions of Proposition 3.2, there exists a positive constant $C>0$ which is independent of $\alpha, \beta$ and $T$  such that for $0\leq t\leq T$,
\begin{equation}\label{3.33}
  \aligned
     &\|(\phi,\psi,\phi_x)(t)\|^2+\int_0^t\left\|(\sqrt V_x\psi,\psi_x)(\tau)\right\|^2\ d\tau \\
     \le&C\left(\|(\phi_0,\psi_0,\phi_{0x})\|^2+e^{-c_1\beta}+(\delta^2+\varepsilon)\int_0^t\|(\phi_x,\phi_{xx},\psi_{xx})(\tau)\|^2\ d\tau\right)
       \endaligned
\end{equation}
holds provided that $\epsilon$ and $\delta$ are suitably small.
\end{Lemma}
\noindent{\bf Proof.}~~Multiplying  $(\ref{3.3})_1$ by $\phi$, $(\ref{3.3})_2$ by $-\frac{1}{p^{\prime}(V)}\psi$, then adding these equalities and integrating the resulting equation over $[0,T]\times\mathbb R^+$, we have
\begin{equation}\label{3.34}
  \begin{split}
     & \int_0^{\infty}\left(\frac{\phi^2}{2}+\frac{\psi^2}{2(-p^{\prime}(V))}\right)\ dx+\int_0^t\int_0^{\infty}\frac{sp^{\prime\prime}(V)V_x\psi^2}{2(p^{\prime}(V))^2}\ dxd\tau+\int_0^t\int_0^{\infty}\frac{\mu \psi_x^2}{-p^{\prime}(V)V}\ dxd\tau \\
      =& \int_0^{\infty}\left(\frac{\phi_0^2}{2}+\frac{\psi_0^2}{2(-p^{\prime}(V_0))}\right)dx+\mu\int_0^t\int_0^{\infty}\left(\frac{1}{p^{\prime}(V)V}\right)_x\psi_x\psi dxd\tau\\
      &+\kappa\int_0^t\int_0^{\infty}\frac{\phi_{xxx}\psi}{p^{\prime}(V)v^5}\ dxd\tau+\mu\int_0^t\int_0^{\infty}\frac{U_x\phi_x\psi}{p^{\prime}(V)V^2}\ dxd\tau+\int_0^t\int_0^{\infty}\frac{-F\psi}{p^{\prime}(V)}\ dxd\tau\\
      &+\mu\int_0^t\frac{\psi_x\psi}{p^{\prime}(V)V}|_{x=0}\ d\tau+\int_0^t\psi\phi|_{x=0} d\tau \\
      =&\int_0^{\infty}\left(\frac{\phi_0^2}{2}+\frac{\psi_0^2}{2(-p^{\prime}(V_0))}\right) dx+\sum_{i=1}^6I_i,
  \end{split}
\end{equation}
where $V_0=V(x+\alpha-\beta)$.
It follows from Proposition 2.1,  (\ref{3.14}),  (\ref{3.15}) and the Cauchy inequality that
\begin{equation}\label{3.35}
  \begin{split}
    |I_1|+|I_3| & \le C\int_0^t\int_0^{\infty}(|V_x\psi\psi_x|+|V_x\phi_x\psi|)dxd\tau \\
      & \le \frac{1}{8}\int_0^t\int_0^{\infty}\frac{sp^{\prime\prime}(V)V_x\psi^2}{(p^{\prime}(V))^2}\ dxd\tau+C\int_0^t\int_0^{\infty}\|V_x(\tau)\|_{L^{\infty}}(\psi_x^2+\phi_x^2)\ dxd\tau\\
      &\le \frac{1}{8}\int_0^t\int_0^{\infty}\frac{sp^{\prime\prime}(V)V_x\psi^2}{(p^{\prime}(V))^2}\ dxd\tau+C\delta^2\int_0^t\int_0^{\infty}\left(\psi_x^2+\phi_x^2\right)\ dxd\tau,
  \end{split}
\end{equation}
\begin{equation}\label{3.36}
  \begin{split}
    I_2 &=\kappa \int_0^t\int_0^{\infty}\left\{\left(\frac{\kappa \phi_{xx}\psi}{p^{\prime}(V)v^5}\right)_x-(\frac{\psi}{p^{\prime}(V)v^5})_x\phi_{xx}\right\}dxd\tau\\
      &= -\kappa \int_0^t\frac{\phi_{xx}\psi}{p^{\prime}(V)v^5}\bigg|_{x=0}\ d\tau-\kappa\int_0^t\int_0^{\infty}\left(\frac{\psi_x\phi_{xx}}{p^{\prime}(V)v^5}+\psi\frac{-(p^{\prime}(V)v^5)_x}{(p^{\prime}(V)v^5)^2}\phi_{xx}\right)dxd\tau,
  \end{split}
\end{equation}
Using $(\ref{3.3})_1$ and integration by parts, we have
\begin{eqnarray}\label{3.37}
    -\kappa\int_0^t\int_0^{\infty}\frac{\psi_x\phi_{xx}}{p^{\prime}(V)v^5}\ dxd\tau
    &&= -\kappa\int_0^t\int_0^{\infty}\frac{\phi_t\phi_{xx}}{p^{\prime}(V)v^5}\ dxd\tau\nonumber \\
      &&= -\kappa\int_0^t\int_0^{\infty}\left[\left(\frac{\phi_t\phi_{x}}{p^{\prime}(V)v^5}\right)_x-\left(\frac{\phi_t}{p^{\prime}(V)v^5}\right)_x\right]\phi_xdxd\tau\nonumber \\
      &&=\kappa\int_0^t\frac{\phi_t\phi_{x}}{p^{\prime}(V)v^5}\bigg|_{x=0}\ d\tau+\kappa\int_0^t\int_0^{\infty}\frac{\phi_{tx}\phi_{x}}{p^{\prime}(V)v^5}\ dxd\tau\nonumber \\
      &&\quad-\kappa\int_0^t\int_0^{\infty}\frac{\psi_x\phi_x(p^{\prime}(V)v^5)_x}{(p^{\prime}(V)v^5)^2}\ dxd\tau\nonumber\\
      &&=\kappa\int_0^t\frac{\phi_t\phi_{x}}{p^{\prime}(V)v^5}\bigg|_{x=0}\ d\tau+\kappa\int_0^{\infty}\frac{\phi_x^2}{2p^{\prime}(V)v^5}\ dx-\kappa\int_0^{\infty}\frac{\phi_{0x}^2}{2p^{\prime}(V_0)v_0^5}dx\nonumber \\
      &&\quad+\frac{\kappa}{2}\int_0^t\int_0^{\infty}\phi_x^2\frac{-sV^{\prime}p^{\prime\prime}(V)}{(p^{\prime}(V))^2v^5}\ dxd\tau+\frac{5\kappa}{2}\int_0^t\int_0^{\infty}\phi_x^2\frac{U^{\prime}+\psi_x}{p^{\prime}(V)v^6}dxd\tau\nonumber \\
      &&\quad-\kappa\int_0^t\int_0^{\infty}\frac{\psi_x\phi_x(p^{\prime}(V)v^5)_x}{(p^{\prime}(V)v^5)^2}\ dxd\tau,
\end{eqnarray}
Combining (\ref{3.36})-(\ref{3.37}) yields
\begin{equation}\label{3.38}
  \begin{split}
    I_2= & \kappa\int_0^{\infty}\frac{\phi_x^2}{2p^{\prime}(V)v^5}\ dx-\kappa\int_0^{\infty}\frac{\phi_{0x}^2}{2p^{\prime}(V_0)v_0^5}\ dx-\kappa\int_0^t\frac{\phi_{xx}\psi}{p^{\prime}(V)v^5}\bigg|_{x=0}d\tau+\kappa\int_0^t\frac{\psi_{x}\phi_x}{p^{\prime}(V)v^5}|_{x=0}d\tau \\
      &+\frac{\kappa}{2}\int_0^t\int_0^{\infty}\left(\phi_x^2\frac{-sV^{\prime}p^{\prime\prime}(V)}{(p^{\prime}(V))^2v^5}+5\phi_x^2\frac{U^{\prime}+\psi_{xx}}{p^{\prime}(V)v^6}\right) dxd\tau\\
      &+\kappa\int_0^t\int_0^{\infty}\frac{p^{\prime\prime}(V)V^{\prime}v^5+p^{\prime}(V)5v^4(V^{\prime}+\phi_{xx})}{(p^{\prime}(V)v^5)^2}\phi_{xx}\psi\ dxd\tau\\
      &-\kappa\int_0^t\int_0^{\infty}\psi_x\phi_x\frac{p^{\prime\prime}(V)V^{\prime}v^5+p^{\prime}(V)5v^4(V^{\prime}+\phi_{xx})}{(p^{\prime}(V)v^5)^2}\ dxd\tau\\
      =&\kappa\int_0^{\infty}\frac{\phi_x^2}{2p^{\prime}(V)v^5}\ dx-\kappa\int_0^{\infty}\frac{\phi_{0x}^2}{2p^{\prime}(V_0)v_0^5}\ dx+\sum_{i=1}^5I_{2i},
  \end{split}
\end{equation}
From (\ref{3.16}), (\ref{3.17}), Proposition 2.1,  (\ref{3.14}) and (\ref{3.15}), we obtain
\begin{equation}\label{3.39}
    |I_{21}|+ |I_{22}|+|I_5|+|I_6| \le Ce^{-c_1\beta},
\end{equation}
Using  Proposition 2.1,  the Sobolev inequality, the  Cauchy inequality,  (\ref{3.14}), (\ref{3.15}) and the a priori assumption (\ref{3.13}), we get
  \begin{equation}\label{3.40}
\begin{split}
      I_{23}&\le C\int_0^t\int_0^{\infty}\left(\|(V^{\prime},U^{\prime})(\tau)\|_{L^{\infty}}\phi_x^2+\|\phi_x(\tau)\|_{L^{\infty}}|\phi_x\psi_{xx}|\right) dxd\tau\\
      &\le C\delta^2\int_0^t\int_0^{\infty}\phi_x^2\ dxd\tau+CN(T)\int_0^t\int_0^{\infty}\left(\phi_x^2+\psi_{xx}^2\right)dxd\tau\\
      &\le C(\delta^2+\varepsilon)\int_0^t\int_0^{\infty}\left(\phi_x^2+\psi_{xx}^2\right) dxd\tau,\\
       \end{split}
 \end{equation}
\begin{equation}\label{3.41}
\begin{split}
      I_{24}&\le \eta\int_0^t\int_0^{\infty}V^{\prime}\psi^2\ dxd\tau+C_{\eta}\int_0^t\int_0^{\infty}(V^{\prime}\phi_{xx}^2+\|\psi(\tau)\|_{L^{\infty}}\phi_{xx}^2) dxd\tau\\
      &\le \eta\int_0^t\int_0^{\infty}V^{\prime}\psi^2\ dxd\tau+C_{\eta}(\delta^2+\varepsilon)\int_0^t\int_0^{\infty}\phi_{xx}^2 dxd\tau.
  \end{split}
\end{equation}
Here and hereafter,  $\eta$ is a small positive constant and $C_\eta$ is a positive constant depending on $\eta$.
Similarly, we have
\begin{equation}\label{3.42}
  I_{25}\le C(\delta^2+\varepsilon)\int_0^t\int_0^{\infty}\left(\psi_x^2+\phi_x^2+\phi_{xx}^2\right) dxd\tau.
\end{equation}
Thus it follows from (\ref{3.38})-(\ref{3.42}) that
\begin{equation}\label{3.43}
   \aligned
  I_2&\le \kappa\int_0^{\infty}\frac{\phi_x^2}{2p^{\prime}(V)v^5}\ dx-\kappa\int_0^{\infty}\frac{\phi_{0x}^2}{2p^{\prime}(V_0)v_0^5}\ dx+\eta\int_0^t\int_0^{\infty}V^{\prime}\psi^2\ dxd\tau+
  \\
 &\quad+C_{\eta}(\delta^2+\varepsilon)\int_0^t\int_0^{\infty}(\psi_x^2+\phi_x^2+\phi_{xx}^2+\psi_{xx}^2)\ dxd\tau+Ce^{-c_1\beta}.
   \endaligned
\end{equation}
Due to the proof of Proposition 2.1 in \cite{Chen-2012-3}, $V^\prime(x-st+\alpha-\beta)>0$ and $\lim_{x\rightarrow\pm\infty}\frac{V^{\prime\prime}}{V^{\prime}}=\lambda_{\pm}$, where $\lambda_{\pm}\in\mathbb R$ are two constants, thus it follows from Proposition 2.1 that
\begin{equation}\label{3.44}
  \begin{split}
    F & =O(1)\left(|\phi_x^2|+|\phi_{xx}^2|+|\psi_{xx}\phi_x|+)\left(\left|\frac{V^{\prime\prime}}{V^{\prime}}\right|+|V^{\prime}|\right)|V^{\prime}\phi_x|+|V^{\prime}\phi_{xx}|\right)\\
      & =O(1)\left(|\phi_x^2|+|\phi_{xx}^2|+|\psi_{xx}\phi_x|+|V^{\prime}\phi_x|+|V^{\prime}\phi_{xx}|\right).
  \end{split}
\end{equation}
Then similar to the estimates of $I_2$, we have
\begin{equation}\label{3.45}
  \begin{split}
    |I_4| & \le C\int_0^t\int_0^{\infty}|F||\psi|\ dxd\tau \\
      & \le C\int_0^t\int_0^{\infty}(|\phi_x^2|+|\psi_{xx}\phi_x|+|\phi_{xx}^2|+|V^{\prime}\phi_x|+|V^{\prime}\phi_{xx}|)|\psi|\ dxd\tau\\
      &\le \frac{1}{8}\int_0^t\int_0^{\infty}\frac{sp^{\prime\prime}(V)V_x\psi^2}{(p^{\prime}(V))^2}\ dxd\tau+C(\delta^2+\varepsilon)\int_0^t\int_0^{\infty}(\phi_x^2+\phi_{xx}^2+\psi_{xx}^2)\ dxd\tau.
        \end{split}
\end{equation}
Inserting (\ref{3.35}),  (\ref{3.39}), (\ref{3.43}) and (\ref{3.45}) into (\ref{3.34}), and using Proposition 2.1, (\ref{3.14})-(\ref{3.15}) and the smallness of $\eta,\varepsilon$ and $\delta$,  we can get (\ref{3.33}). Thus the proof of Lemma 3.3 is completed.

The next lemma give the estimate on $\|\phi_x(t)\|$.
\begin{Lemma}
 Under the assumptions of Proposition 3.2, there exists a positive constant $C>0$ which is independent of $\alpha, \beta$ and $T$  such that for $0\leq t\leq T$, it holds that
\begin{equation}\label{3.46}
\|\phi_x(t)\|^2+\int_0^t\|\phi_x(\tau)\|_1^2\ d\tau\le C\left(\|(\phi_{0x},\psi_0,\phi_0)\|^2+e^{-c_1\beta}+(\delta^2+\varepsilon)\int_0^t\|\psi_{xx}(\tau)\|^2 d\tau\right)
\end{equation}
 provided that $\epsilon$ and $\delta$ are suitably small.
\end{Lemma}
\noindent{\bf Proof.}~~Multiplying $(\ref{3.3})_2$ by $-\phi_x$, and integrating the resulting equation with respect to  $x$ over $[0,+\infty)$ yields
\begin{equation}\label{3.47}
  \begin{split}
     &\int_0^{\infty}-p^{\prime}(V)\phi_x^2\ dx\underbrace{-\kappa\int_0^{\infty}\frac{\phi_{xxx}\phi_x}{v^5}\ dx}_{I_7} \\
     = & \underbrace{\int_0^{\infty}\psi_t\phi_x\ dx}_{I_8}  \underbrace{-\int_0^{\infty}\frac{\mu}{V}\psi_{xx}\phi_x\ dx}_{I_9}  \underbrace{-\mu\int_0^{\infty}\frac{U^{\prime}\phi_x^2}{V^2}\ dx}_{I_{10}}   \underbrace{-\int_0^{\infty}F\phi_x\ dx}_{I_{11}}.
  \end{split}
\end{equation}
By Proposition 2.1, (\ref{3.14})-(\ref{3.15}) and integrations by parts, we have
\begin{eqnarray}
     I_7 &&=-\kappa\int_0^{\infty}\left(\phi_{xx}\frac{\phi_x}{v^5}\right)_x\ dx+\kappa\int_0^{\infty}\frac{\phi_{xx}^2}{v^5}\ dx-5\kappa\int_0^{\infty}\frac{\phi_{xx}(\phi_x+V_x)\phi_x}{v^6} dx\nonumber \\
       && \ge\kappa\phi_{xx}\frac{\phi_x}{v^5}\bigg|_{x=0}+\kappa\int_0^{\infty}\frac{\phi_{xx}^2}{v^5}\ dx-C(\delta^2+\varepsilon)\int_0^{\infty}(\phi_{xx}^2+\phi_x^2)dx, \\
        I_8 &&= \int_0^{\infty}[(\psi\phi_x)_t-\psi\phi_{xt}]dx\nonumber \\
      &&=\frac{d}{dt}\int_0^{\infty}\psi\phi_x\ dx-\int_0^{\infty}[(\psi\psi_x)_x-\psi_x^2] dx\nonumber \\
      &&=\frac{d}{dt}\int_0^{\infty}\psi\phi_x\ dx+\psi\psi_x|_{x=0}+\int_0^{\infty}\psi_x^2 dx,\\
      I_9&&=-\mu\int_0^{\infty}\frac{\phi_{tx}\phi_x}{V}\ dx=-\frac{\mu}{2}\frac{d}{dt}\int_0^{\infty}\frac{\phi_x^2}{V}\ dx+\frac{\mu}{2}\int_0^{\infty}\frac{V^{\prime}}{V^2}\phi_x^2dx\nonumber \\
      &&\le -\frac{\mu}{2}\frac{d}{dt}\int_0^{\infty}\frac{\phi_x^2}{V}\ dx+C\delta^2\int_0^{\infty}\phi_x^2\ dx,\\
      I_{10} &&\le C\delta^2\int_0^{\infty}\phi_x^2\ dx,\\
      I_{11} && =\int_0^{\infty}|F||\phi_x|\ dx\le O(1)\int_0^{\infty}(|\phi_x^2|+|\psi_{xx}\phi_x|+|\phi_{xx}^2|+|V^{\prime}\phi_x|+|V^{\prime}\phi_{xx}|)|\phi_x|dx\nonumber \\
      && \le O(1)(\delta^2+\varepsilon)\int_0^{\infty}(\phi_x^2+\psi_{xx}^2+\phi_{xx}^2)dx.
\end{eqnarray}
Substituting  (3.50)-(3.54) into (\ref{3.47}), we get by the smallness of $\delta$ and $\varepsilon$ that
\begin{equation}\label{3.53}
  \begin{split}
     & \frac{\mu}{2}\frac{d}{dt}\int_0^{\infty}\frac{\phi_x^2}{V}\ dx+\int_0^{\infty}-p^{\prime}(V)\phi_x^2\ dx+\kappa\int_0^{\infty}\frac{\phi_{xx}^2}{v^5}\ dx \\
      \le & \frac{d}{dt}\int_0^{\infty}\psi\phi_x\ dx+\psi\psi_x|_{x=0}-\kappa\frac{\phi_x\phi_{xx}}{v^5}\bigg|_{x=0}+\int_0^{\infty}\psi_x^2dx.
  \end{split}
\end{equation}
Integrating (\ref{3.53}) in $t$ over $[0,+\infty)$ gives
\begin{equation}\label{3.54}
  \begin{split}
     & \frac{\mu}{4}\int_0^{\infty}\frac{\phi_x^2}{V}\ dx+\int_0^t\int_0^{\infty}\left(-p^{\prime}(V)\phi_x^2+\kappa\frac{\phi_{xx}^2}{v^5}\right)dxd\tau \\
     \le & C(\|(\phi_{0x},\psi_0)\|^2+\|\psi(t)\|^2)+\int_0^t\|\psi_x(\tau)\|^2\ d\tau+\int_0^t\psi\psi_x|_{x=0}\ d\tau-\kappa\int_0^t\frac{\phi_x\phi_{xx}}{v^5}\bigg|_{x=0}\ d\tau,
  \end{split}
\end{equation}
where we have used the fact that
\[
  \int_0^{\infty}\psi\phi_xdx\le \frac{\mu}{4}\int_0^{\infty}\frac{\phi_x^2}{V}\ dx+C\|\psi(t)\|^2.
\]

Thus (\ref{3.46}) follows from (\ref{3.54}), (\ref{3.16}), (\ref{3.17}), Lemma 3.3 and the smallness of $\varepsilon$ and $\delta$ immediately.  This completes the proof of Lemma 3.4.

Next, we estimate $\|\psi_x(t)\|$.
\begin{Lemma}
 Under the assumptions of Proposition 3.2, there exists a positive constant $C>0$ which is independent of $\alpha, \beta$ and $T$  such that for $0\leq t\leq T$, it holds that
\begin{equation}\label{3.57}
\aligned
  &\|(\psi_x,\phi_{xx})(t)\|^2+\int_0^t\|\psi_{xx}(\tau)\|^2\ d\tau\\
  &\le C\left(\|\phi_0\|_2^2+\|\psi_0\|_1^2+e^{-c_1\beta}+(\delta^2+\varepsilon)\int_0^t\|\phi_{xxx}(\tau)\|^2\ d\tau\right)
\endaligned
\end{equation}
 provided that $\epsilon$ and $\delta$ are suitably small.
\end{Lemma}
\noindent{\bf Proof.}~~Multiplying $(\ref{3.3})_2$ by $-\psi_{xx}$, and integrating the resulting equation in  $x$ over $[0,+\infty)$, we have
\begin{equation}\label{3.58}
 \aligned
    \frac{1}{2}\frac{d}{dt}\int_0^{\infty}\psi_x^2\ dx+\int_0^{\infty}\frac{\mu}{V}\psi_{xx}^2 dx
   & = \frac{1}{2}\int_0^{\infty}\psi_{0x}^2\ dx+\underbrace{\int_0^{\infty}p^{\prime}(V)\phi_x\psi_{xx}\ dx}_{J_1} -\psi_t\psi_x|_{x=0}\\
      &\quad +\underbrace{\kappa\int_0^{\infty}\frac{\phi_{xxx}\psi_{xx}}{v^5}\ dx}_{J_2}+\underbrace{\mu\int_0^{\infty}\frac{U^{\prime}\phi_x\psi_{xx}}{V^2}\ dx}_{J_3}\underbrace{-\int_0^{\infty}F\psi_{xx}\ dx}_{J_4},
\endaligned
\end{equation}
The Cauchy inequality implies that
\begin{equation}\label{3.59}
    J_1  \le \eta\int_0^{\infty}\psi_{xx}^2\ dx+C_{\eta}\int_0^{\infty}(p^{\prime}(V))^2\phi_x^2dx \le\eta\int_0^{\infty}\psi_{xx}^2\ dx+C_{\eta}\int_0^{\infty}\phi_x^2\ dx,
\end{equation}
\begin{equation}\label{3.60}
      J_3 \le \eta\int_0^{\infty}\psi_{xx}^2\ dx+C_{\eta}\delta^2\int_0^{\infty}\phi_x^2dx.
\end{equation}
Using $(\ref{3.3})_1$, we get
\[
  \begin{split}
      J_2& =\kappa\int_0^{\infty}\frac{\phi_{xxx}\phi_{tx}}{v^5}\ dx=\kappa\int_0^{\infty}[(\phi_{xx}\frac{\phi_{tx}}{v^5})_x-\phi_{xx}(\frac{\phi_{tx}}{v^5})_x]\ dx\\
      &=-\frac{\kappa}{2}\frac{d}{dt}\int_0^{\infty}\frac{\phi_{xx}^2}{v^5}\ dx-\kappa\phi_{xx}\frac{\psi_{xx}}{v^5}\bigg|_{x=0}\underbrace{-\frac{5\kappa}{2}\int_0^{\infty}\phi_{xx}^2\frac{\psi_{xx}+U^{\prime}}{v^6}\ dx}_{J_{21}}\\
      &\quad+\underbrace{5\kappa\int_0^{\infty}\phi_{xx}\psi_{xx}\frac{\phi_{xx}+V^{\prime}}{v^6}\ dx}_{J_{22}}.
  \end{split}
\]
It follows from  Proposition 2.1,  the Sobolev inequality, the  Cauchy inequality,  (\ref{3.14}), (\ref{3.15}) and the a priori assumption (\ref{3.13}) that
\begin{equation*}
  \begin{split}
          J_{21} &\le C\int_0^{\infty}\|\phi_{xx}(t)\|_{L^{\infty}}|\phi_{xx}\psi_{xx}|+\|U^{\prime}(t)\|_{L^{\infty}}\phi_{xx}^2\ dx \\
            &\le C\left(\delta^2\|\phi_{xx}(t)\|^2+\|\phi_{xx}(t)\|^{\frac{3}{2}}\|\phi_{xxx}(t)\|^{\frac{1}{2}}\|\psi_{xx}(t)\|\right),\\
            &\le C(\delta^2+\varepsilon)\|(\phi_{xx},\phi_{xxx},\psi_{xx})(t)\|^2,
        \end{split}
\end{equation*}

\begin{equation*}
  J_{22}\le C(\delta^2+\varepsilon)\|(\phi_{xx},\phi_{xxx},\psi_{xx})(t)\|^2.
\end{equation*}
Thus
\begin{equation}\label{3.61}
    J_2 \le -\frac{\kappa}{2}\frac{d}{dt}\int_0^{\infty}\frac{\phi_{xx}^2}{v^5}\ dx-\kappa\phi_{xx}\frac{\psi_{xx}}{v^5}|_{x=0}+C(\delta^2+\varepsilon)\|(\phi_{xx},\phi_{xxx},\psi_{xx})(t)\|^2. \\
\end{equation}
Using  (\ref{3.44}), $J_4$ can be controlled as follows:
\begin{eqnarray}\label{3.62}
     J_4 &&\le \int_0^{\infty}|F\psi_{xx}| dx\nonumber\\
     &&\le \eta\int_0^{\infty}\psi_{xx}^2\ dx+C_{\eta}\int_0^{\infty}|F|^2dx\nonumber\\
     &&\le \eta\int_0^{\infty}\psi_{xx}^2\ dx+C_{\eta}\int_0^{\infty}\left(\|\phi_{x}(t)\|_{L^{\infty}}^2(\phi_x^2+\psi_{xx}^2)+\|V^{\prime}(t)\|_{L^{\infty}}^2(\phi_x^2+\phi_{xx}^2)+\|\phi_{xx}(t)\|_{L^{\infty}}^2\phi_{xx}^2\right)dx\nonumber\\
     &&\le \eta\|\psi_{xx}(t)\|^2+C_{\eta}(\delta^2+\varepsilon)\|(\phi_x,\psi_{xx},\phi_{xx})(t)\|^2.
\end{eqnarray}
Inserting (\ref{3.59})-(\ref{3.62}) into (\ref{3.58}), and integrating the resulting formula in $t$ over $[0,t]$, we have by the smallness of $\eta,\delta$ and $\varepsilon$ that
\begin{equation}\label{3.63}
  \begin{split}
     & \|\psi_{x}(t)\|^2+\|\phi_{xx}(t)\|^2+\int_0^t\|\psi_{xx}(\tau)\|^2\ d\tau \\
     \le & C\left(\|(\psi_{0x},\phi_{0xx})\|^2+\int_0^t\|\phi_x(\tau)\|_1^2\ d\tau\right)+C(\delta^2+\varepsilon)\int_0^t\|\phi_{xxx}(\tau)\|^2\ d\tau\\
     &+\int_0^t\left(-\psi_t\psi_x-\kappa\phi_{xx}\frac{\psi_{xx}}{v^5}\right)\bigg|_{x=0}\ d\tau.
  \end{split}
\end{equation}
Then (\ref{3.57}) can be obtained by (\ref{3.63}), (\ref{3.16}),(\ref{3.17}), Lemma 3.4  and the smallness of $\varepsilon$ and $\delta$. This completes the proof of Lemma 3.5.

For the estimates on $\displaystyle\int_0^t\|\phi_{xx}(\tau)\|_1^2\ d\tau$, we have
\begin{Lemma}
 Under the assumptions of Proposition 3.2, there exists a positive constant $C>0$ which is independent of $\alpha, \beta$ and $T$  such that for $0\leq t\leq T$,
\begin{equation}\label{3.64}
  \|\phi_{xx}(t)\|^2+\int_0^t\|\phi_{xx}(\tau)\|_1^2 d\tau\le C\left(\|\phi_0\|_2^2+\|\psi_0\|_1^2+e^{-c_1\beta}\right)
\end{equation}
 holds provided that $\epsilon$ and $\delta$ are suitably small.
\end{Lemma}
\noindent{\bf Proof.}~~Multiplying $(\ref{3.3})_2$ by $\phi_{xxx}$ and utilizing $(\ref{3.3})_1$, we have
\begin{equation}\label{3.65}
  \begin{split}
     & -p^{\prime}(V)\phi_{xx}^2+\kappa\frac{\phi_{xxx}^2}{v^5}+\left(\frac{\mu\phi_{xx}^2}{2V}\right)_t-(\psi_x\phi_{xx})_t \\
     = &\left (\frac{\mu}{V}\psi_{xx}\phi_{xx}\right)_x-\left(\frac{\mu}{V}\right)_x\phi_{xx}\psi_{xx}+\left(\frac{\mu}{V}\right)_t\frac{\phi_{xx}^2}{2}-(p^{\prime}(V)\phi_x\phi_{xx})_x+p^{\prime\prime}(V)V_x\phi_x\phi_{xx}\\
     & -(\psi_t\phi_{xx})_x-(\psi_x\psi_{xx})_x+\psi_{xx}^2-\mu\frac{U^{\prime}\phi_x}{V^2}\phi_{xxx}+F\phi_{xxx}
      \end{split}
\end{equation}
Integrating (\ref{3.65}) over $[0,t]\times\mathbb R$ and using  the Cauchy inequality, we obtain
\begin{eqnarray}\label{3.66}
     && \frac{1}{4}\int_{\mathbb R}\frac{\mu\phi_{xx}^2}{V}\ dx+\int_0^t\int_{\mathbb R}\left(-p^{\prime}(V)\phi_{xx}^2+\frac{\kappa\phi_{xxx}^2}{v^5}\right) dxd\tau\nonumber \\
     \le && C\|(\phi_{0xx},\psi_{0x})\|^2+C\|\psi_x(t)\|^2+\left|\int_0^tK_1\ d\tau\right|+\int_0^t\|\psi_{xx}(\tau)\|^2\ d\tau\nonumber \\
     && +\left|\int_0^t\psi_t\phi_{xx}|_{x=0}\ d\tau\right|+\int_0^t\int_{\mathbb R}\left(|K_2|+|K_3|\right) dxd\tau,
\end{eqnarray}
where
\begin{equation*}
  \begin{split}
    K_1 =& \left(-\frac{\mu}{V}\psi_{xx}\phi_{xx}+\psi_x\psi_{xx}+p^{\prime}(V)\phi_x\phi_{xx}\right)\bigg|_{x=0}, \\
    K_2 = &\frac{\mu V_x}{V^2}\phi_{xx}\psi_{xx}-\frac{\mu}{2V^2}V_t\phi_{xx}^2+p^{\prime\prime}(V)V_x\phi_x\phi_{xx},\\
    K_3=&F\phi_{xxx}-\mu\frac{U^{\prime}\phi_x\phi_{xxx}}{V^2}.
  \end{split}
\end{equation*}
Similar to the estimates of $J_{21}$,  we have
\begin{equation}\label{3.67}
  \begin{split}
    \int_0^t\int_{\mathbb R}|K_2|\ dxd\tau \le & C\int_0^t\int_0^{\infty}\left(|V_x\phi_{xx}\psi_{xx}|+|V^{\prime}\phi_{xx}^2|+|V_x\phi_x\phi_{xx}|\right)\ dxd\tau \\
     \le & C\delta^2\int_0^t\|(\phi_{xx},\psi_{xx},\phi_x)(\tau)\|^2\ d\tau,
  \end{split}
\end{equation}
\begin{equation}\label{3.68}
  \begin{split}
     \int_0^t\int_{\mathbb R}|K_3|\ dxd\tau \le & C\int_0^t\int_{\mathbb R}\left(|\phi_x^2|+|\psi_{xx}\phi_x|+|\phi_{xx}^2|+|V^{\prime}(\phi_x,\phi_{xx})|\right)|\phi_{xxx}|\ dxd\tau \\
     \le & C(\delta^2+\varepsilon)\int_0^t\|(\phi_x,\psi_{xx},\phi_{xx},\phi_{xxx})(\tau)\|^2\ d\tau.
  \end{split}
\end{equation}
Combining (\ref{3.66})-(\ref{3.68}),   using (\ref{3.16})-(\ref{3.17}) and the smallness of $\delta,\varepsilon$ yields  (\ref{3.64}) at once. This finishes the proof of Lemma 3.6.

As a consequence of Lemmas 3.3-3.6,  we have
\begin{Corollary}
Under the assumptions of Proposition 3.2, there exists a positive constant $C>0$ which is independent of $\alpha, \beta$ and $T$  such that
\begin{equation}\label{3.69}
  \begin{split}
     & \|\phi(t)\|_2^2+\|\psi(t)\|_1^2+\int_0^t\left(\|\phi_{x}(\tau)\|_2^2+\|\psi_x(\tau)\|_1^2\right) d\tau \\
     \le & C\left(\|\phi_0\|_2^2+\|\psi_0\|_1^2+e^{-c_1\beta}\right), \quad \forall\, t\in[0,T].
  \end{split}
\end{equation}
\end{Corollary}

For the estimates on $\displaystyle\|\psi_{xx}(\tau)\|$, we have
\begin{Lemma}
 Under the assumptions of Proposition 3.2, there exists a positive constant $C>0$ which is independent of $\alpha, \beta$ and $T$  such that for $0\leq t\leq T$,
\begin{equation}\label{3.70}
  \begin{split}
     & \|(\phi_{xxx},\psi_{xx})(t)\|^2+\int_0^t\|\psi_{xxx}(\tau)\|^2\ d\tau \\
     \le & C\left(\|\phi_0\|_3^2+\|\psi_0\|_2^2+e^{-c_1\beta}+\varepsilon^{\frac{1}{2}}\int_0^t\|\phi_{xxx}(\tau)\|^2\ d\tau\right)
  \end{split}
\end{equation}
 holds provided that $\epsilon$ and $\delta$ are suitably small.
\end{Lemma}
\noindent{\bf Proof.}~~Differentiating  $(\ref{3.3})_2$  with respect to $x$ once,  then multiplying  the resulting equation  by $-\psi_{xxx}$ and integrating over $[0,t]\times\mathbb R$, we get
\begin{eqnarray}\label{3.71}
     &&\frac{1}{2}\int_0^{\infty}\psi_{xx}^2\ dx+\int_0^t\int_0^{\infty}\frac{\mu}{V}\psi_{xxx}^2\ dxd\tau\nonumber \\
    =&&\frac{1}{2}\int_0^{\infty}\psi_{0xx}^2\ dx-\int_0^t\psi_{xt}\psi_{xx}|_{x=0}\ d\tau+\kappa\int_0^t\int_0^{\infty}\left(\frac{\phi_{xxx}}{v^5}\right)_x\psi_{xxx}dxd\tau-\int_0^t\int_0^{\infty}F_x\psi_{xxx}\ dxd\tau\nonumber \\
    &&+\int_0^t\int_0^{\infty}Q_1dxd\tau,
\end{eqnarray}
where\[Q_1=p^{\prime\prime}(V)V_x\phi_x\psi_{xxx}-p^{\prime}(V)\phi_{xx}\psi_{xxx}-\left(\frac{\mu}{V}\right)_x\psi_{xx}\psi_{xxx}-\mu\left(\frac{U^{\prime}\phi_x}{v^2}\right)_x\psi_{xxx}.
\]
Using integration by parts, we have
\begin{equation}\label{3.72}
  \begin{split}
     &\kappa\int_0^t\int_0^{\infty}\left(\frac{\phi_{xxx}}{v^5}\right)_x\psi_{xxx}\ dxd\tau \\
    =&\kappa\int_0^t\int_0^{\infty}\left(\frac{\phi_{xxx}}{v^5}\right)_x\phi_{txx}\ dxd\tau\\
    =&\kappa\int_0^t\int_0^{\infty}\left[\left(\frac{\phi_{xxx}}{v^5}\phi_{txx}\right)_x-\frac{\phi_{xxx}}{v^5}\phi_{txxx}\right]dxd\tau\\
    =&-\kappa\int_0^t\frac{\phi_{xxx}}{v^5}\phi_{txx}\bigg|_{x=0}\ d\tau-\kappa\int_0^t\int_0^{\infty}\left[\left(\frac{\phi_{xxx}^2}{2v^5}\right)_t+\frac{5\phi_{xxx}^2v_t}{2v^6}\right]dxd\tau   \\
    =&-\kappa\int_0^t\frac{\phi_{xxx}}{v^5}\phi_{txx}\bigg|_{x=0}\ d\tau-\kappa\int_0^{\infty}\frac{\phi_{xxx}^2}{2v^5}\ dx
    +\kappa\int_0^{\infty}\frac{\phi_{0xxx}^2}{2v_0^5}\ dx-\frac{5\kappa}{2}\int_0^t\int_0^{\infty}\frac{\phi_{xxx}^2u_x}{v^6}\ dxd\tau.
  \end{split}
\end{equation}
It follows from the Sobolev inequality and Proposition 2.1 that
\begin{equation}\label{3.73}
  \begin{split}
    \left|\int_0^t\int_0^{\infty}\frac{\phi_{xxx}^2u_x}{v^6}\ dxd\tau\right| \le &C\int_0^t\int_0^{\infty}\left|\phi_{xxx}^2(U^{\prime}+\psi_{xx})\right| dxd\tau  \\
     \le & C\delta^2\int_0^t\|\phi_{xxx}(\tau)\|^2\ d\tau+C\int_0^t\|\psi_{xx}(\tau)\|^{\frac{1}{2}}\|\psi_{xxx}(\tau)\|^{\frac{1}{2}}\|\phi_{xxx}(\tau)\|^2\ d\tau\\
     \le& C\delta^2\int_0^t\|\phi_{xxx}(\tau)\|^2\ d\tau+\eta\int_0^t\|\psi_{xxx}(\tau)\|^2\ d\tau\\
     &+C_{\eta}\int_0^t\left(\|\psi_{xx}(\tau)\|^2+\|\phi_{xxx}(\tau)\|^4\right)d\tau.
  \end{split}
\end{equation}
On the other hand, since
\begin{equation*}
  \begin{split}
    F_x= & -\{p(v)-p(V)-p^{\prime}(V)\phi_x\}_x-\left(\frac{\mu\psi_{xx}\phi_x}{vV}\right)_x \\&
    +\left\{\kappa V^{\prime\prime}\left(\frac{1}{V^5}-\frac{1}{v^5}\right) +\frac{5\kappa}{2}\frac{\phi_{xx}^2+2\phi_{xx}V^{\prime}}{v^6}+\frac{5\kappa}{2}(V^{\prime})^2\left(\frac{1}{v^6}-\frac{1}{V^6}\right)+\frac{\mu U^{\prime}\phi_x^2}{vV^2}\right\}_x\\
      =&O(1)\left(|(V_x,U^{\prime\prime})|\phi_x^2+|\phi_{xx}\phi_x^2|+|\phi_x\phi_{xx}|+|\psi_{xxx}\phi_x|+|\psi_{xx}\phi_{xx}|+|V_x\phi_x\psi_{xx}|\right.\\
      &+|\phi_{xx}\phi_x\psi_{xx}|+|(V^{\prime\prime\prime},V^{\prime\prime}V^{\prime},(V^{\prime})^3)||\phi_x|+|(V^{\prime\prime},(V^{\prime})^2)||\phi_{xx}|+|V^{\prime}\phi_{xxx}|\\
      &\left.+|\phi_{xx}\phi_{xxx}|+|\phi_{xx}^3|+|\phi_{xx}^2V^{\prime}|\right),
  \end{split}
\end{equation*}
we obtain  from the Cauchy inequality, the Sobolev inequality,  Proposition 2.1 and the a priori assumption (\ref{3.13}) that
\begin{eqnarray}\label{3.74}
     &&\left|- \int_0^t\int_0^{\infty}F_x\psi_{xxx}\ dxd\tau\right|\nonumber \\
     \le && \eta\int_0^t\int_0^{\infty}\psi_{xxx}^2\ dxd\tau+C_{\eta}\int_0^t\int_0^{\infty}|F_x|^2\ dxd\tau\nonumber \\
     \le &&\eta\int_0^t\int_0^{\infty}\psi_{xxx}^2\ dxd\tau+C_{\eta}\int_0^t\int_0^{\infty}\left\{|(V^{\prime},U^{\prime\prime})|^2\phi_x^4+|\phi_{xx}|^2\phi_x^4+|\phi_x^2\phi_{xx}^2|+|\psi_{xx}^2\phi_x^2|+|\psi_{xx}^2\phi_{xx}^2|\right.\nonumber \\
     &&+|V_x^2\phi_x^2\psi_{xx}^2|+|\phi_{xx}\phi_x\psi_{xx}|^2+|(V^{\prime\prime\prime},V^{\prime\prime}V^{\prime},(V^{\prime})^3)|^2|\phi_x^2|+|(V^{\prime\prime},(V^{\prime})^2)|^2|\phi_{xx}|^2+|V^{\prime}\phi_{xxx}|^2\nonumber \\
      &&\left.+|\phi_{xx}\phi_{xxx}|^2+|\phi_{xx}^6|+|\phi_{xx}^4V^{\prime 2}|\right\}dxd\tau\nonumber \\
      \le &&\eta\int_0^t\|\psi_{xxx}(\tau)\|^2d\tau+C_{\eta}(\delta^2+\varepsilon^2)\int_0^t\|(\phi_x, \phi_{xx}, \phi_{xxx}, \psi_{xxx})(\tau)\|^2d\tau.
\end{eqnarray}
Similarly, it holds
\begin{equation}\label{3.75}
      \int_0^t\int_0^{\infty}|Q_1| dxd\tau
     \le  \eta\int_0^t\|\psi_{xxx}(\tau)\|^2 d\tau+C_{\eta}\int_0^t\|(\phi_x,\phi_{xx}, \psi_{xx})(\tau)\|^2d\tau.
\end{equation}
Putting (\ref{3.72})-(\ref{3.75}) into (\ref{3.71}), then (\ref{3.70}) follows  from Corollary 3.1, (\ref{3.16}), (\ref{3.18}) and the smallness of $\delta,\varepsilon$ and $\eta$.
This completes the proof of Lemma 3.7.

Finally, we give the estimates of $\displaystyle\int_0^t\|\phi_{xxx}(\tau)\|_1^2\ d\tau$.
\begin{Lemma}
 Under the assumptions of Proposition 3.2, there exists a positive constant $C>0$ which is independent of $\alpha, \beta$ and $T$  such that for $0\leq t\leq T$, it holds
\begin{equation}\label{3.76}
  \|\phi_{xxx}(t)\|^2+\int_0^t\|\phi_{xxx}(\tau)\|_1^2\ d\tau\le C\left(\|\phi_0\|_3^2+\|\psi_0\|_2^2+e^{-c_1\beta}\right)
\end{equation}
 provided that $\epsilon$ and $\delta$ are suitably small.
\end{Lemma}
\noindent{\bf Proof.}~~Differentiating  $(\ref{3.3})_2$  with respect to $x$ once,  then multiplying  the resulting equation  by $\phi_{xxx}$ and integrating over $[0,t]\times\mathbb R$
 gives
\begin{equation}\label{3.77}
  \begin{split}
     & \frac{1}{2}\int_0^{\infty}\frac{\mu}{V}\phi_{xxx}^2dx+\int_0^t\int_0^{\infty}\left(-p^{\prime}(V)\phi_{xxx}^2+\frac{\kappa\phi_{xxxx}^2}{v^5}\right)\,dxd\tau-\int_0^{\infty}\psi_{xx}\phi_{xxx}\,dx \\
     = & \frac{1}{2}\int_0^{\infty}\frac{\mu}{V_0}\phi_{0xxx}^2dx-\int_0^{\infty}\psi_{0xx}\phi_{0xxx}\,dx+\int_0^t\int_0^{\infty}\sum_{i=2}^4Q_i\,dxd\tau,
  \end{split}
\end{equation}
where
\begin{equation*}
  \begin{split}
    Q_2= & \left(p^{\prime}(V)\phi_{xx}\phi_{xxx}+\psi_{tx}\phi_{xxx}+\psi_{xx}\psi_{xxx}-\frac{\mu}{V}\phi_{txx}\phi_{xxx}\right)\bigg|_{x=0}, \\
    Q_3=  & p^{\prime\prime}(V)V_x\phi_x\phi_{xxx}+p^{\prime}(V)\phi_{xxx}^2- \left(\frac{\mu}{V}\right)_x\psi_{xxx}\phi_{xxx}+ \left(\frac{\mu}{V}\right)_t\frac{\phi_{xxx}^2}{2},\\
    Q_4=&-p^{\prime}(V)\phi_{xx}\phi_{xxx}-p^{\prime\prime}(V)V^{\prime}\phi_x\phi_{xxxx}+\left(\frac{\mu}{V}\right)_x\psi_{xx}\phi_{xxxx}+5\kappa\phi_{xxx}\phi_{xxxx}\frac{v_x}{v^6}\\
    &-\mu\left(\frac{U^{\prime}\phi_x}{V^2}\right)_x\phi_{xxxx}+F_x\phi_{xxxx}.
  \end{split}
\end{equation*}
The estimates (\ref{3.18})-(\ref{3.18-3}) imply that
\begin{equation}\label{3.78}
  \begin{split}
    \left|\int_0^t\int_0^{\infty}Q_2dxd\tau\right|\le &2\eta\int_0^t\|\phi_{xxxx}(\tau)\|^2\ d\tau+C_{\eta}\left(\int_0^t\|\phi_{xx}(\tau)\|_1^2\ d\tau+e^{-c_1\beta}\right) \\
      & +C\left(\int_0^t\|\psi_{xx}(\tau)\|_1^2\ d\tau+\varepsilon^{\frac{1}{2}}\int_0^t\|\phi_{xxxx}(\tau)\|^2\ d\tau\right).
  \end{split}
\end{equation}
Similar to (\ref{3.73})-(\ref{3.74}), we have
\begin{equation}\label{3.79}
   \left|\int_0^t\int_0^{\infty}Q_3 dxd\tau\right|\le C\int_0^t\|(\phi_x,\phi_{xxx},\psi_{xxx})(\tau)\|^2\ d\tau,
\end{equation}
\begin{equation}\label{3.80}
   \left|\int_0^t\int_0^{\infty}Q_4dxd\tau\right|\le\eta\int_0^t\int_0^{\infty}\phi_{xxxx}^2\ dxd\tau+C_{\eta}\int_0^t(\|\phi_x(\tau)\|_2^2+\|\psi_{xx}(\tau)\|^2)\ d\tau,
\end{equation}
and
\begin{equation}\label{3.81}
  \left|-\int_0^{\infty}\psi_{xx}\phi_{xxx} dx\right|\le \frac{1}{4}\int_0^{\infty}\frac{\mu}{V}\phi_{xxx}^2\ dx+C\|\psi_{xx}(t)\|^2.
\end{equation}
Combining (\ref{3.77})-(\ref{3.81}),  and using Corollary 2.1, Lemma 3.7 and   the smallness of $\varepsilon$ and $\eta$, we obtain (\ref{3.76}).
This completes the proof of Lemma 3.8.

\noindent{\bf Proof of Proposition 3.2.}~~Proposition 3.2 follows  from Corollary 2.1 and Lemmas 3.7-3.8 immediately.

\subsection{Proof of Theorem 2.1}
To complete the proof  Theorem 2.1,  we first observe that  the a priori assumption (\ref{3.13}) can be closed by choosing $\|\phi_0\|, \|\psi_0\|_1$  and $\beta^{-1}$  sufficiently  small such that
$$\|\phi_0\|^2_2+\|\psi_0\|^2_1+e^{-c_1\beta}<\frac{\epsilon^2}{4C_0}.$$
Then based on Propositions 3.1-3.2,  the standard continuation argument asserts that there exists a unique  global (in time) solution $(\phi,\psi)(t,x)\in X_{\tilde{M}}(0,+\infty)$
to the initial-boundary value  problem (\ref{3.4})-(\ref{3.9}), where $\hat{M}=2\sqrt{C_0(\|\phi_0\|_3^2+\|\psi_0\|_2^2+e^{-c_1\beta})}$. Moreover, we can derive from  (\ref{3.11}) and  the system (\ref{3.3}) that
\begin{equation}\label{3.82}
\int_0^{+\infty}\left(\|\phi_x(t)\|_2^2+\|\psi_x(t)\|_1^2+\left|\frac{d}{dt}\left(\|\phi_x(t)\|_2^2+\|\psi_x(t)\|_1^2\right) \right|\right)\,dt<\infty,
\end{equation}
which implies that $\|\phi_x(t)\|_2+ \|\psi_x(t)\|_1\rightarrow0$ as $t\rightarrow+\infty$.  Then by the Sobolev inequality, we have
\begin{equation}\label{3.83}
\|\phi_x(t)\|_{L^\infty}\leq\|\phi_x(t)\|^{\frac{1}{2}}\|\phi_{xx}(t)\|^{\frac{1}{2}}\rightarrow0,\mbox{ as}\,\, t\rightarrow+\infty.
 \end{equation}
and
 \begin{equation}\label{3.84}
\|\psi_x(t)\|_{L^\infty}\leq\|\psi_x(t)\|^{\frac{1}{2}}\|\psi_{xx}(t)\|^{\frac{1}{2}}\rightarrow0,\mbox{ as}\,\, t\rightarrow+\infty.
 \end{equation}
Thus (\ref{3.12}) is proved and the proof of Theorem 2.1 is completed.

\bigbreak

\begin{center}
{\bf Acknowledgement}
\end{center}
The first author would like to thank professor Huijiang Zhao for some useful discussions on  the topic of this paper.  
This work  was supported by the National Natural Science Foundation of China (Grant No. 11501003),   the  Doctoral Scientific Research Fund of Anhui University (Grant No. J10113190005), and the Cultivation Fund of  Young Key Teacher at Anhui University.

\end{document}